\documentclass[10pt,leqno]{amsart}
\usepackage{amssymb,amsbsy,amsmath,amsfonts,amssymb,amscd,times,graphics,color,xypic}

\newtheorem{theorem}{Theorem}

\newtheorem{mainproposition}{Main Proposition}
\newtheorem{lemma}{Lemma}
\newtheorem{mainlemma}{Main Lemma}
\newtheorem{corollary}{Corollary}
\newtheorem{definition}{Definition}

\usepackage{graphics,epsfig,amssymb,amsfonts,amscd}
\def\C{{\mathbb C}}

\def\K{{\mathbb K}}

\def\N{{\mathbb N}}

\def\R{{\mathbb R}}
\def\V{{\mathbb V}}

\begin{document}


\title[Propagation of analyticity for essentially finite smooth 
CR mappings]{
Propagation of analyticity for essentially 
finite $\mathcal{ C}^\infty$-smooth 
CR mappings
}

\author{Jo\"el Merker}

\address{
CNRS, Universit\'e de Provence, LATP, UMR 6632, CMI, 
39 rue Joliot-Curie, F-13453 Marseille Cedex 13, France}

\email{merker@cmi.univ-mrs.fr} 

\subjclass[2000]{Primary: 32H02. Secondary: 32H04, 32V20, 32V30, 
32V35, 32V40}

\date{\number\year-\number\month-\number\day}

\begin{abstract}
An analytico-geometric reflection principle is established by means of
normal deformations of analytic discs.
\end{abstract}

\maketitle

\begin{center}
\begin{minipage}[t]{11cm}
\baselineskip =0.35cm
{\scriptsize

\centerline{\bf Table of contents}

\smallskip

{\bf 1.~Introduction \dotfill 1.}

{\bf 2.~Preliminaries and statement of the main theorem \dotfill 1.}

{\bf 3.~Polynomial identities \dotfill 4.}

{\bf 4.~Proof of Lemma~2.11 \dotfill 5.}

{\bf 5.~Description of the proof of Theorem~2.9 \dotfill 6.}

{\bf 6.~Rings of CR functions and their fields of fractions 
\dotfill 8.}

{\bf 7.~Statement of Main Lemma 7.1 \dotfill 11.}

{\bf 8.~Localization at a nice boundary point \dotfill 14.}

{\bf 9.~Construction of analytic discs \dotfill 16.}

{\bf 10.~Meromorphic extension \dotfill 17.}

{\bf 11.~Proof of Lemma~2.12 \dotfill 20.}

{\bf 12.~References \dotfill 21.}

\smallskip

{\footnotesize\tt [With 5 figures]}

}\end{minipage}
\end{center}

\bigskip

\section*{\S1.~Introduction}

The analyticity of local $\mathcal{ C}^\infty$-smooth CR
diffeomorphisms between two essentially finite generic real analytic
submanifolds of $\C^n$ is established in~\cite{ bjt}, as a kind
of reflection principle, provided that all the components of the CR
diffeomorphism extend holomorphically to a fixed wedge 
(the so-called notion of {\sl essential finiteness} appeared
in the work~\cite{ dw} by K. Diederich and S.M.~Webster;
a further geometric reflection principle for locally finite
$\mathcal{ C}^\infty$-smooth CR mappings appeared in 
the work~\cite{ df2} by K.~Diederich and J.E.~Forn{\ae}ss).
In \cite{
br1}, \cite{ br2} and more recently in~\cite{ cps1}, \cite{ cps2},
\cite{ da}, \cite{ mmz1} ({\it cf.}~also the applications \cite{
cdms}, \cite{ mmz2}), separate assumptions on the map
and on the target have been
unified: instead, it is assumed that the 
so-called {\sl characteristic
variety} is zero-dimensional at the central point. However, 
in these references, it is
always supposed that the source generic submanifold is minimal at the
central point, whereas, in~\cite{ bjt}, no minimality assumption was
needed. This article is devoted to fill the gap between these two
trends of thought, applying the technique of normal deformations of
analytic discs borrowed from~\cite{ tu2},
\cite{ mp1}, \cite{ mp2}.

\section*{\S2.~Preliminaries and statement of the main theorem}

\subsection*{2.1.~Initial data} 
Let $\K = \R$ or $\C$. Let $\nu \in \N$. Let ${\sf x} \in \K^\nu$.
Set $\vert {\sf x} 
\vert := \max_{1\leq i\leq \nu} \vert {\sf x}_i 
\vert$. For
$\rho >0$, denote $\Delta_n (\rho):= \{ {\sf x } \in \K^\nu : \ \vert
{\sf x} \vert < \rho\}$.

Consider a local 
$\mathcal{ C}^\infty$-smooth CR mapping between
two local generic submanifolds $M$ in $\C^n$ and $M'$ in $\C^{n'}$, 
defined {\it precisely}\, as follows (background
material may be found in~\cite{ ber}, \cite{ m2}).
The purpose is to avoid the (ambiguous) language of germs.

\def\thedefinition{2.2}\begin{definition}
{\rm 
A {\sl local $\mathcal{ C}^\infty$-smooth 
CR mapping} consists of the following data.
\begin{itemize}
\item[{\bf (I)}]
A {\sl local generic submanifold $M$ in $\C^n$} of positive
codimension $d\geq 1$ and of positive CR dimension $m:=n-d \geq 1$
passing through a point $p_0\in \C^n$ and defined in coordinates
$t=(z,w) =(x+iy,u +iv)\in \C^m\times \C^d$ vanishing at $p_0$ as a
graph:
\def\theequation{2.3}\begin{equation}
M=\{(z,w)\in \Delta_n(\rho_1): \, 
v_j=\varphi_j(x,y,u), \, j=1,\dots,d\},
\end{equation}
where the functions $\varphi_j$ are real analytic on $\Delta_{ 2m+d}
(2\rho_1)$, for some $\rho_1>0$. It is
also required for all $\rho$
with $0\leq \rho \leq 2\rho_1$, that $\vert \varphi (x,y,u) \vert <
\rho$ if $\vert (x,y,u) \vert < \rho$, namely $M$ is a uniformly
approximatively horizontal graph. Of course, after perhaps shrinking
$\rho_1>0$, this condition is automatically satisfied if the
coordinates are adjusted 
at the beginning in order that $T_0M=\{v=0\}$. In fact, it
is more convenient to work with a local representation of $M$ by
complex defining equations
\def\theequation{2.4}\begin{equation}
M=\{(z,w)\in\Delta_n(\rho_1): \, 
\bar w_j=\Theta_j(\bar z, z, w), \, 
j=1,\dots,d\},
\end{equation}
obtained by applying the implicit function theorem to~\thetag{ 2.3},
where of course, one
may assume that the $\Theta_j$ converge normally for
$\vert (\bar z, z, w)\vert < 2\rho_1$ and that for all $\rho$ with
$0\leq \rho \leq 2\rho_1$, one has $\vert \Theta_j(\bar z, z, w) \vert
< \rho$ if $\vert (\bar z,z,w) \vert < \rho$.
\item[{\bf (II)}]
A local generic submanifold $M'$ in $\C^{n'}$ with central point
$p_0'$ of positive codimension $d'\geq 1$ and of positive CR dimension
$m':= n' -d' \geq 1$ passing through a point $p_0' \in \C^n$ and
defined in coordinates $t' = (z', w') \in \C^{m'} \times \C^{ d'}$ 
similarly as in~\thetag{ 2.4}.
\item[{\bf (III)}] 
A $\mathcal{ C }^\infty$-smooth mapping $t' = h(t)= (f(t),g(t))=
(z',w')$ with $h(p_0) =p_0'$ which is defined in $M\cap \Delta_n
(\rho_1)\equiv M$ and which satisfies for some two radii $\rho_2$,
$\rho_2'$ with $0< \rho_2 < \rho_1$, $0 < \rho_2' < \rho_1'$, the
condition
\def\theequation{2.5}\begin{equation}
h\left( M\cap \Delta_n(\rho_2) \right)
\subset M'\cap
\Delta_{n'}'(\rho_2').
\end{equation}
\end{itemize}
}
\end{definition}

By~\cite{ bt} (and also \cite{ ber}), after 
shrinking (if necessary) $\rho_1>0$ and
$\rho_2>0$ with $0< \rho_2< \rho_1$:
\begin{itemize}
\item[{\bf (IV)}] 
every $\mathcal{ C}^\infty$-smooth CR function defined on $M\cap
\Delta_n(\rho_1)$ (and in particular the $n'$ components $h_1,\,
\dots,\, h_{n'}$ of $h$) is a uniform limit of polynomials on $M \cap
\Delta_n(\rho_2)$.
\end{itemize}

\def\thedefinition{2.6}\begin{definition}
{\rm 
A {\sl complete wedge in $\Delta_n (\rho_2)$ with edge $M \cap
\Delta_n (\rho_2)$} is a subset of $\C^n$ of the form $\mathcal{ W} =
\mathcal{ W}(M, C, \Delta_n (\rho_2)) = \{(z,w) \in \Delta_n(\rho_2);
\, v-\varphi( x,y,u)\in C\}$, where $C$ is some open strictly convex
infinite ({\it i.e.}~not truncated) cone in $\R^d$.
}
\end{definition}

As in~\cite{ bjt}, it will be assumed that:

\begin{itemize}
\item[{\bf (V)}]
there exists a complete wedge $\mathcal{ W}_2$ in $\Delta_n( \rho_2)$
with edge $M \cap \Delta_n( \rho_2)$ such that the $n$ components of
$h$ extend holomorphically to $\mathcal{ W}_2$.
\end{itemize}

\subsection*{2.7.~CR differentiations}
Put $r_j (t, \bar t):= \bar w_j - \Theta_j (\bar z, z, w)$ for $j=1,
\dots, d$ and and $r_{ j'}' (t', \bar t') := \bar w_{j '}' - \Theta_{
j'} ' ( \bar z', z', w')$ for $j'= 1,\dots, d'$.  Let $\overline{ L
}_1, \dots, \overline{ L}_m$ be an arbitrary basis of $(0,1)$ vector
fields tangent to $M$ having real analytic coefficients (the most
convenient is written in~\thetag{ 3.2} below). Consider the
{\sl first characteristic variety of $h$ at $p_0$} to be the complex
analytic subset $\V_0'$ of $\Delta_{ n'}'( \rho_2')$ consisting of
elements $t'$ satisfying the equations
\def\theequation{2.8}\begin{equation}
\left.
\left[ \overline{L}^\beta 
r_{j'}'(t', \overline{h(t)})\right]
\right\vert_{\bar t=0}=0, \, \ 
\text{\rm for all} \ j'=1,\dots,d' \
\text{\rm and all} \ \beta\in\N^m.
\end{equation}
It is indeed a complex analytic subset defined as the zero set of an
infinite collection of functions which are holomorphic in $\Delta_{n'}
(\rho_2')$. By~\thetag{ 2.5}, $r_{j'}' \left( h(t),\overline{ h(t)}
\right) =0$, for $j'=1,\dots,d'$ and for all $t \in M$. It follows
that the origin $p_0'$ belongs to the complex analytic subset
$\V_0'$. The focus is on the dimension at $p_0 '$ of $\V_0 '$. The map
$h$ will be called {\sl essentially 
finite at} $p_0$ if $\dim_{p_0 '} \V_0 '
=0$.  Denote by $\mathcal{ O}_{ CR} (M, p_0)$ the (not local) CR orbit
of $p_0$ in $M$. The main result is as follows.

\def\thetheorem{2.9}\begin{theorem}
Let $h: M \to M'$ be a local $\mathcal{ C }^\infty$-smooth CR mapping
between two real analytic local generic submanifolds of $\C^n$ and of
$\C^{ n'}$. Assume that there exists a complete wedge $\mathcal{ W}_2$
in $\Delta_n( \rho_2)$ with edge $M\cap \Delta_n( \rho_2)$ such that
the $n'$ components of $h$ extend holomorphically to $\mathcal{ W}_2$,
assume that there exist points
\def\theequation{2.10}\begin{equation}
q_0\in {\mathcal{O}}_{CR}(M,p_0)
\cap \Delta_n(\rho_2)
\end{equation}
arbitrarily close to $p_0$ at which $h$ is real analytic and assume
that $h$ is essentially finite at $p_0$. Then there exists a radius
$\rho_3>0$ with $0< \rho_3 < \rho_2 < \rho_1$ such that $h$ extends
holomorphically to $\Delta_n(\rho_3)$.
\end{theorem}

In the version of Theorem~2.9 published in~\cite{ da}, \cite{ mmz1},
it is assumed that $M$ is minimal at $p_0$, which entails the
holomorphic extendability assumption {\bf (V)}, thanks to~\cite{
tu1}. In~\cite{ da}, \cite{ mmz1}, there is
a crucial proposition about
envelopes of meromorphy (same statement and same proof in the two
references\,--\,albeit notations
differ), relying on subtle geometric
arguments which stem from the theory of 
deformations of analytic discs
developed in~\cite{ a}, \cite{ tu1}, 
\cite{ mp1}.  In this
article, the purpose is to clean up and to simplify these geometric
arguments, by means of the propagation of wedge extendability theorem
established in~\cite{ tr}, \cite{ tu2} ({\it cf.}~\cite{ ht} for a
preliminary version). The stronger Theorem~2.9 will be established
thanks to this change of geometric point of view.  An elementary lemma
applies to recover from Theorem~2.9 
the main result of~\cite{ da}, \cite{ mmz1}.

\def\thelemma{2.11}\begin{lemma}
Let $h: M \to M'$ be a local $\mathcal{ C}^\infty$-smooth CR mapping
between two real analytic local generic submanifolds of $\C^n$ and of
$\C^{ n'}$. Assume that $M$ is minimal at $p_0$, so that $\mathcal{
O}_{ CR} (M, p_0)$ contains $M \cap \Delta_n (\rho_2)$ for some
$\rho_2 >0$ and so that {\rm (}thanks to~{\rm \cite{ tu1})} after
perhaps shrinking $\rho_2>0$, the assumption {\bf (V)} holds. If $h$
is essentially finite at $p_0$, then there exist points $q_0 
\in M \cap \Delta_n (\rho_2)$ arbitrarily close
to $p_0$ at which $h$ is real analytic.
\end{lemma}

Finally, in order to recover the main result of~\cite{ bjt}, remind
that the essential finiteness of $M'$ at $p_0'$ together with the CR
diffeomorphism assumption entails the essential 
finiteness of $h$ at $p_0$
(\cite{ da}, Lemma~4.1; a more general version is Corollary~1.3
in~\cite{mmz1}; the most general version appears as Theorem~4.3.1 {\bf
(3)} in~\cite{ m2}, in which it is shown that CR-transversality of $h$
at $p_0$ together with essential finiteness of $M'$ at $p_0'$ implies
that $h$ is essentially finite at $p_0$). In~\cite{ bjt}, it is
observed that essential finiteness of a hypersurface $M$ at one of its
points $p_0$ implies its minimality (finite type in the sense of
Lie-Chow-Kohn-Bloom-Graham) at $p_0$.
Here are further observations.

\def\thelemma{2.12}\begin{lemma}
If $M'$ is a local generic submanifold of $\C^{ n'}$ passing through
a point $p_0'$ which is essentially finite at $p_0'$, then $\dim_\R
\, \mathcal{ O}_{ CR} (M', p_0') \geq 2 {\rm CRdim} \, M' + 1$ and the
CR orbit $\mathcal{ O}_{ CR} (M', p_0')$ itself is essentially finite
at $p_0'$. Furthermore, in the case where $M$ is a real analytic
hypersurface, essential finiteness of $h$ at $p_0$ implies that $M$ is
minimal at $p_0$.
\end{lemma}

Assume that $n=n'$ and that $h$ is a CR diffeomorphism.  Then there
exists a Zariski-dense open subset of points $q_0' \in \mathcal{ O}_{
CR} (M', p_0')$ at which $M'$ is finitely nondegenerate.  It follows
that $h$ itself is finitely nondegenerate at $q_0 := h^{ -1} (q_0')$,
and by a known result (\cite{ p1}, \cite{ ha}, \cite{
ber}, \cite{ la}), $h$ is real analytic at $q_0$.
Applying then Theorem~2.9: 

\def\thecorollary{2.13}\begin{corollary}
{\rm (\cite{ bjt})} Let $h: M \to M'$ be a local $\mathcal{
C}^\infty$-smooth CR mapping which is a CR diffeomorphism. If
$M$ is essentially finite at $p_0$ and if the
components of $h$ extend holomorphically to a wedge at $p_0$, then $h$
is real analytic at $p_0$.
\end{corollary}

Further applications (in the spirit of~\cite{ cdms}, \cite{ mmz2})
that may be stated are left to the interested reader. The remainder of
this article is devoted to the proofs of Theorem~2.9, of Lemma~2.11
and of Lemma~2.12.

\section*{ \S3.~Polynomial identities}

\subsection*{ 3.1.~Differentiations}
Denote by $\overline{L}_1,\dots,\overline{L}_m$ the basis of
$(0,1)$ vector fields tangent to $M$ defined explicitely by
\def\theequation{3.2}\begin{equation}
\overline{L}_k:=
{\partial \over \partial \bar z_k}+
\sum_{j=1}^d\, {\partial \Theta_j \over
\partial \bar z_k} \left(\bar z,z,w \right)\, 
{\partial \over \partial \bar w_j}.
\end{equation}
Here, the coefficients of the vector fields $\overline{L}_k$ are
holomorphic with respect to $w$.  Since it will be more convenient for
later use to have antiholomorphic dependence with respect to $w$,
replace $w$ by $\overline{\Theta}(z,\bar z,\bar w)$ (which is possible
when $(z,w)$ belongs to $M$), and write the vector fields under the
form
\def\theequation{3.3}\begin{equation}
\overline{L}_k:=
{\partial \over \partial \bar z_k}+
\sum_{j=1}^d\, {\partial \Theta_j \over
\partial \bar z_k}
\left(\bar z,z,
\overline{\Theta}(z,\bar z,\bar w)
\right)\, 
{\partial \over \partial \bar w_j}.
\end{equation}
Apply the derivations $( \overline{ L})^\beta$, $\beta \in \N^m$,
to $r_{j'}' \left( h(t), \overline{ h(t)} \right) =0$, which yields
\def\theequation{3.4}\begin{equation}
(\overline{L })^\beta \, r_{j'}' \left( h(t), 
\overline{h (t)} \right) =0,
\end{equation}
for $t\in M \cap \Delta_n ( \rho_2)$ and for $j'=1, \dots,d'$. 

As the coefficients of the vector fields $\overline{L}_k$ are
holomorphic with respect to $(z,\bar t)$, the differentiated
equations~\thetag{3.4} may be rewritten under the
developed form:
\def\theequation{3.5}\begin{equation}
R_{j',\beta}' \left(z,\bar t, 
J_{\bar t}^{\vert \beta \vert}
\overline{h(t)} :\, h(t) \right) =0.
\end{equation}
Here, $J_{ \bar t}^\ell \overline{ h} (\bar t) : = \left(
\partial_{\bar t}^\alpha \overline{ h_{ i'}( t)} \right)_{ 1\leq
i'\leq n', \, \vert \alpha \vert \leq \ell }$ denotes the $\ell$-th
jet of $\overline{ h} (\bar t)$. By construction, the functions
$R_{ j', \beta}' = R_{ j', \beta}' \left( z, \bar t, \overline{ J}^{
\vert \beta \vert} :\, t' \right)$ are holomorphic with respect to
$z$, $\bar t$ with $\vert z \vert < \rho_2$, $\vert \bar t \vert <
\rho_2$, they are holomorphic with respect to the zero-th order jet
$J_{\bar t}^0 \overline{ h}( \bar t) \equiv \overline{ h}(\bar t)$
with $\vert J_{ \bar t}^0 \vert < \rho_2'$, they are relatively
polynomial with respect to the nonzero derivatives $\left( \partial_{
\bar t}^\alpha \overline{ h_{ i'}( t)} \right)_{ 1 \leq i'\leq n', \,
1 \leq \vert \alpha \vert \leq \vert \beta \vert}$, and they are
relatively holomorphic with respect to $t'$ with $\vert t' \vert <
\rho_2'$.

By the main assumption of essential finiteness, there exists an integer
$\ell_0\geq 1$ such that the complex analytic subset
defined by the equations
\def\theequation{3.6}\begin{equation}
R_{j',\beta}'
\left(
0,0,J_{\bar t}^{\vert \beta \vert} \overline{h}(0) :\,
t'\right)
=0, \ \ \ \ \ j'=1,\dots,d', \ \vert \beta \vert \leq \ell_0,
\end{equation}
which passes through the origin in $\C^{ n'}$, 
is of dimension zero at the origin.

By~\cite{ ber}, chapter~5, it follows that there exists $\rho_3$ with
$0< \rho_3 < \rho_2$, there exists $\rho_3'$ with $0 < \rho_3' <
\rho_2'$, there exists $\varepsilon$ with $\varepsilon >0$ and for all
$i' = 1, \dots, n'$, there exist monic Weierstrass polynomials $P_{
i'}\left( z, \bar t, J_{ \bar t}^{ \ell_0 }: t_{ i'}' \right)$ of the form
\def\theequation{3.7}\begin{equation}
P_{i'} \left( z,\bar t, 
J_{\bar t}^{\ell_0} :\, t_{i'}' \right) =
(t_{i'}')^{N_{i'}'}+
\sum_{1\leq I' \leq N_{i'}'}\, 
H_{i',I'} \left( z,\bar t, 
J_{\bar t}^{\ell_0} \right)\, 
(t_{i'}')^{N_{i'}'-I'},
\end{equation} 
with coefficients $H_{i',I'}$ being holomorphic with respect to $z$,
$\bar t$ with $\vert z \vert < \rho_3$, $\vert \bar t \vert < \rho_3$
and with respect to $J_{\bar t}^{\ell_0}$ with $\left\vert J_{\bar
t}^{\ell_0} - J_{\bar t}^{\ell_0} \overline{ h}(0) \right\vert <
\varepsilon$, with moreover $\left\vert J_{\bar t}^{\ell_0}
\overline{h(t)} - J_{ \bar t}^{\ell_0} \overline{h}(0) \right\vert <
\varepsilon$ for all $t \in M\cap \Delta_n(\rho_3)$, such that the
complex analytic set
\def\theequation{3.8}\begin{equation}
\left\{ \left( z,\bar t, J_{\bar t}^{\ell_0},t' \right): \,
R_{j',\beta}'\left(z,\bar t, J_{\bar t}^{\vert \beta \vert} :\, t'
\right)=0, \ j'=1,\dots,d', \ \vert \beta \vert \leq \ell_0\right\}
\end{equation}
is contained in the zero-set of all Weierstrass polynomials 
$P_{i'}$, namely the set:
\def\theequation{3.9}\begin{equation}
\left\{
\left(
z,\bar t, J_{\bar t}^{\ell_0},t' \right) : \, P_{i'} \left( z,\bar t,
J_{\bar t}^{\ell_0} :\, t_{i'}' \right) =0, \ i'=1,\dots,n'
\right\}.
\end{equation}
Thanks to Hilbert's Nullstellensatz, there exists an integer $\nu \geq
1$ such that for $i'=1,\dots,n'$, the powers $(P_{i' })^\nu$ belong to
the ideal generated by the $R_{j', \beta}'$. Consequently, each
component $h_{i'}(t)$ satisfies the monic polynomial equation
\def\theequation{3.10}\begin{equation}
(h_{i'}(t))^{N_{i'}'}+
\sum_{1\leq I' \leq N_{i'}'}\,
H_{i',I'} \left( z,\bar t, 
J_{\bar t}^{\ell_0} \overline{h(t)} \right)\, 
(h_{i'}(t))^{N_{i'}'-I'}= 0,
\end{equation}
for all $t\in M \cap \Delta_n (\rho_3)$. Here, notably, each component
$h_{i'} (t)$ is separated from the others.

\section*{\S4.~Proof of Lemma~2.11} 

To establish Lemma~2.11, let $\rho_2$ be as in its statement and let
$\rho_3$ be as in \S3 just above.  Shrinking $\rho_3$ if necessary,
assume $0 < \rho_3 < \rho_2$ to fix ideas.

Fix $i':= 1$ and consider the following trivial dichotomy: either
\def\theequation{4.1}\begin{equation}
{\partial P_1 \over \partial t_1'} \left( z,\bar t, J_{\bar
t}^{ \ell_0} \overline{ h(t)} :\, h_1(t) \right) = 0,
\end{equation} 
for all $t\in M\cap \Delta_n(\rho_3)$ or
\def\theequation{4.2}\begin{equation}
{\partial P_1 \over \partial t_1'} \left( z,\bar t, J_{\bar
t}^{\ell_0} \overline{h(t)} :\, h_1(t) \right) \not \equiv 0,
\end{equation}
over $M\cap \Delta_n (\rho_3)$. In the first case, replace the
equation $P_1 \left( z, \bar t, J_{\bar t}^{ \ell_0} \overline{ h(t)}
:\, h(t) \right) =0$ by the equation~\thetag{4.1}, which is a
monic
polynomial of degree $N_{i'}'-1$ in $h_{i'}(t)$. After a finite number
of steps, the second case~\thetag{4.2} holds, with a
monic polynomial of
lower degree, still denoted by $P_1$. Pick a point $q_1 \in M \cap
\Delta_n (\rho_3)$ together with an open neighborhood $\omega_1
\subset M \cap \Delta_n (\rho_3)$ of $q_1$ so that ${\partial P_1
\over \partial t_{i'}'} \left( z_{q_1}, \bar t_{q_1}, 
J_{\bar t}^{\ell_0}
\overline{ h(t_{q_1})} :\, h(t_{q_1}) \right) \neq 0$ for all $t \in
\omega_1$, whereas $P_1 \left( z_{q_1}, \bar t_{q_1}, 
J_{\bar t}^{\ell_0}
\overline{ h(t_{q_1})} :\, h(t_{q_1}) 
\right) \equiv 0$ over $\omega_1$.
Fix now $i' := 2$. Apply the same dichotomy as
~\thetag{ 4.1}, \thetag{ 4.2} to $\frac{ 
\partial P_2}{\partial t_2'}$, but for $t$ running only over
$\omega_1$. Process as in the case $i_1'=1$
to replace $P_2$ by a monic polynomial of minimal 
degree, chose $q_2$ and $\omega_2$, {\it etc.}

After $n'$ steps, there exists a point $q_0 \in M\cap 
\Delta_n (\rho_3)$, there exists 
an open neighborhood $\omega_0 \subset
M \cap \Delta_n (\rho_3)$ of $q_0$, there
exists monic polynomials, still denoted
by $P_1,\dots, P_{n '}$ which are
of the form~\thetag{ 3.7}, such that
\def\theequation{4.3}\begin{equation}
P_{i'} \left( z,\bar t, J_{\bar
t}^{ \ell_0} \overline{ h(t)} :\, h_{i'}(t) \right) = 0,
\end{equation} 
for all $t\in \omega_0$ but
\def\theequation{4.4}\begin{equation}
\frac{\partial P_{i'}}{\partial t_{i'}'} \left( z,\bar t, J_{\bar
t}^{\ell_0} \overline{h(t)} :\, h_{i'}(t) \right) \neq 0,
\end{equation}
for all $t \in \omega_0$. 
Apply then the implicit function theorem to the
identities~\thetag{ 4.3}, to solve
\def\theequation{4.5}\begin{equation}
h(t)=\Psi \left( z,\bar t, J_{\bar t}^{
\ell_0} \overline{h(t)} \right),
\end{equation}
for all $t$ in a (possibly smaller) neighborhood of $q_0$, where
$\Psi$ is a certain complex analytic $\C^{n'}$-valued mapping defined
in a neighborhood of $\left( z_{ q_0}, \bar t_{ q_0}, J_{ \bar t}^{
\ell_0} \overline{ h(t_{q_0})} \right)$ in $\C^m \times \C^n \times
\C^{N_{n', n,\ell_0}}$. It then follows from standard arguments which
are easy modifications of the original phenomenon discovered in~\cite{
p1} that $h$ extends holomorphically to a neighborhood of $q_0$ in
$\C^n$, or equivalently that $h$ is real analytic in a neighborhood of
$q_0$ in $M$.

The positive radius $\rho_3$ could have been shrunk from the 
beginning to be arbitrarily small, and the same
reasoning provides a point $q_0$, arbitrarily 
close to $p_0$, at which $h$ is real analytic.

The proof of Lemma~2.11 is complete.
\qed

\section*{ \S5.~Description of the proof of Theorem~2.9}

\subsection*{5.1.~Summary and statement of Main Proposition~5.2}
Let $\rho_1$ be as in Definition~2.2, let $\rho_2$ be as in the
assumptions of Theorem~2.9 and let $\rho_3$ be as in \S3 above.  By
hypothesis, there exists at least one point $q_0 \in \mathcal{ O}_{
CR} (M,p_0) \cap \Delta_n (\rho_3)$ at which $h$ is real analytic (it
will not be necessary to deal with such points which are closer to
$p_0$).

By assumption, the components $h_{i'} (t)$ extend holomorphically to
the complete wedge $\mathcal{ W}_2$ in $\Delta_n (\rho_2)$ with edge
$M\cap \Delta_n (\rho_2)$. Define a complete wedge $\mathcal{ W}_3$
in $\Delta_n (\rho_3)$ with edge $M\cap \Delta_n (\rho_3)$ as the
intersection of $\mathcal{ W}_2$ with $\Delta_n (\rho_3)$.
Here is an illustration, where $\mathcal{ W}_3$
is on the top and its symmetric $\widetilde{ W}_3$ (to be
introduced later) is on the bottom:

\bigskip
\begin{center}
\begin{picture}(0,0)%
\includegraphics{polydiscs-wedges.pstex}%
\end{picture}%
\setlength{\unitlength}{3947sp}%
\begingroup\makeatletter\ifx\SetFigFont\undefined%
\gdef\SetFigFont#1#2#3#4#5{%
  \reset@font\fontsize{#1}{#2pt}%
  \fontfamily{#3}\fontseries{#4}\fontshape{#5}%
  \selectfont}%
\fi\endgroup%
\begin{picture}(5799,3624)(1189,-3973)
\put(1532,-2861){\makebox(0,0)[lb]{\smash{{\SetFigFont{9}{10.8}{\familydefault}{\mddefault}{\updefault}{\color[rgb]{0,0,0}$M$}%
}}}}
\put(6413,-2845){\makebox(0,0)[lb]{\smash{{\SetFigFont{9}{10.8}{\familydefault}{\mddefault}{\updefault}{\color[rgb]{0,0,0}$M$}%
}}}}
\put(6344,-2233){\makebox(0,0)[lb]{\smash{{\SetFigFont{9}{10.8}{\familydefault}{\mddefault}{\updefault}{\color[rgb]{0,0,0}$z,\bar z, {\rm Re} \, w$}%
}}}}
\put(4079,-592){\makebox(0,0)[lb]{\smash{{\SetFigFont{9}{10.8}{\familydefault}{\mddefault}{\updefault}{\color[rgb]{0,0,0}${\rm Im}\, w$}%
}}}}
\put(4049,-2261){\makebox(0,0)[lb]{\smash{{\SetFigFont{9}{10.8}{\familydefault}{\mddefault}{\updefault}{\color[rgb]{0,0,0}$p_0$}%
}}}}
\put(2265,-1158){\makebox(0,0)[lb]{\smash{{\SetFigFont{9}{10.8}{\familydefault}{\mddefault}{\updefault}{\color[rgb]{0,0,0}$\Delta_n(\rho_2)$}%
}}}}
\put(1752,-1038){\makebox(0,0)[lb]{\smash{{\SetFigFont{9}{10.8}{\familydefault}{\mddefault}{\updefault}{\color[rgb]{0,0,0}$\Delta_n(\rho_1)$}%
}}}}
\put(5370,-1165){\makebox(0,0)[lb]{\smash{{\SetFigFont{9}{10.8}{\familydefault}{\mddefault}{\updefault}{\color[rgb]{0,0,0}$\Delta_n(\rho_2)$}%
}}}}
\put(2937,-1453){\makebox(0,0)[lb]{\smash{{\SetFigFont{8}{9.6}{\familydefault}{\mddefault}{\updefault}{\color[rgb]{0,0,0}$\mathcal{W}_3$}%
}}}}
\put(3006,-2770){\makebox(0,0)[lb]{\smash{{\SetFigFont{8}{9.6}{\familydefault}{\mddefault}{\updefault}{\color[rgb]{0,0,0}$\widetilde{\mathcal{W}}_3$}%
}}}}
\put(3450,-1188){\makebox(0,0)[lb]{\smash{{\SetFigFont{9}{10.8}{\familydefault}{\mddefault}{\updefault}{\color[rgb]{0,0,0}$\Delta_n(\rho_3)$}%
}}}}
\put(6004,-1015){\makebox(0,0)[lb]{\smash{{\SetFigFont{9}{10.8}{\familydefault}{\mddefault}{\updefault}{\color[rgb]{0,0,0}$\Delta_n(\rho_1)$}%
}}}}
\put(3484,-2279){\makebox(0,0)[lb]{\smash{{\SetFigFont{9}{10.8}{\familydefault}{\mddefault}{\updefault}{\color[rgb]{0,0,0}$q_0$}%
}}}}
\put(3736,-2300){\makebox(0,0)[lb]{\smash{{\SetFigFont{9}{10.8}{\familydefault}{\mddefault}{\updefault}{\color[rgb]{0,0,0}$\mathcal{V}_4$}%
}}}}
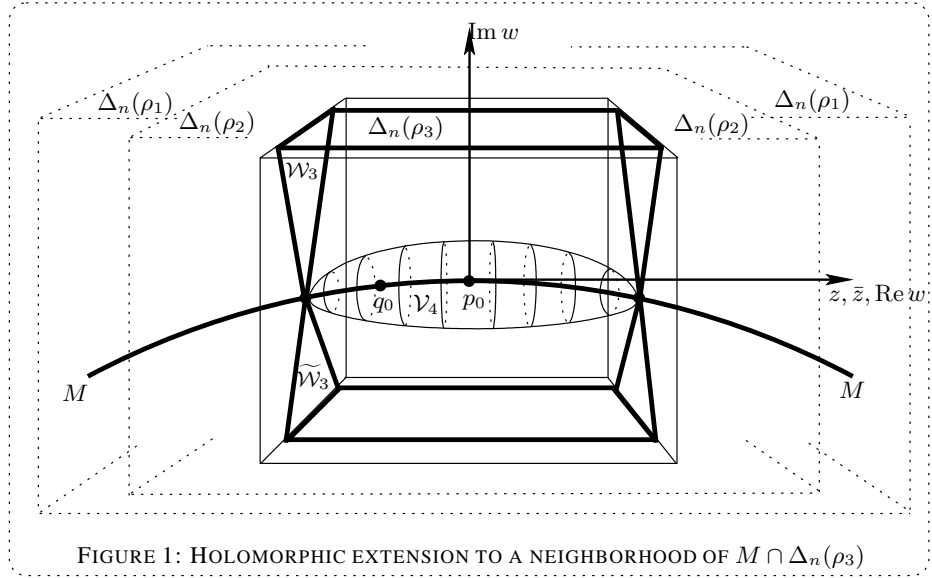
\put(1625,-3881){\makebox(0,0)[lb]{\smash{{\SetFigFont{9}{10.8}{\familydefault}{\mddefault}{\updefault}{\color[rgb]{0,0,0}{\sc Figure~1: Holomorphic extension to a neighborhood of $M\cap \Delta_n(\rho_3)$}}%
}}}}
\end{picture}%

\end{center}
\bigskip

In these conditions, the components $h_{i'} (t)$ of $h$ extend
holomorphically to $\mathcal{W}_3$. Also, the jet $J_{\bar t}^{
\ell_0} \overline{ h(t)}$ in the arguments of the functions $H_{ i',
I'}$ in~\thetag{ 3.10} extends antiholomorphically to $\mathcal{
W}_3$.

It is now time to state an independent general proposition, where $M'$
disapears, where each component $h_{ i'}$ is replaced by a
$\mathcal{ C}^\infty$-smooth CR function
$a$, not necessarily coming from a CR mapping, where the jet $J_{\bar
t}^{ \ell_0} \overline{ h(t)}$ is replaced by an independent vector
valued mapping $\overline{ b}$ which extends antiholomorphically to
$\mathcal{ W}_3$, For similar statements, {\it see}\, \cite{ p3},
\cite{ pu}, \cite{ cps1}, \cite{ cps2}, \cite{ da}, \cite{ mmz1},
\cite{ cdms}, \cite{ mmz2}.

\def\themainproposition{5.2}\begin{mainproposition}
As above, let $M$ be a real analytic local 
generic submanifold defined
as a graph in $\Delta_n (\rho_1)$ by~\thetag{ 2.4}, let $\rho_3$ with
$0< \rho_3 < \rho_1$ and let $\mathcal{ W}_3$ be a complete wedge in
$\Delta_n(\rho_3)$ with edge $M\cap \Delta_n (\rho_3)$. Let $a(t)$ be
a $\mathcal{ C}^\infty$-smooth CR function defined over $M \cap
\Delta_n (\rho_3)$ which extends holomorphically to the complete wedge
$\mathcal{W}_3$ and which is real analytic in a neighborhood of at
least one point $q_0 \in \mathcal{ O}_{ CR}(M,p_0)\cap \Delta_n
(\rho_3)$.  Let $\nu \in\N$, let $\varepsilon >0$ and let $b(t)$ be a
$\C^\nu$-valued $\mathcal{ C}^\infty$-smooth CR mapping defined on $M
\cap \Delta_n(\rho_3)$ which satisfies $\vert b(t) - b(0) \vert <
\varepsilon$ for all $t\in M\cap \Delta_n(\rho_3)$, which extends
holomorphically to $\mathcal{W}_3$ and which is real analytic in a
neighborhood of the same point $q_0\in \mathcal{O}_{CR}(M,p_0)\cap
\Delta_n(\rho_3)$. Let $N \in \N$ with $N\geq 1$, and for $\ell= 0,1,
\dots, N$, let $H_\ell= H_\ell \left( z, \bar t, \bar b\right)$ be
some functions which are holomorphic for $\vert z \vert < \rho_3$, for
$\vert \bar t \vert < \rho_3$ and for $\vert \bar b - \overline{ b(0)}
\vert< \varepsilon$. Assume that $a(t)$ satisfies the {\rm (}not
necessarily monic{\rm )} polynomial equation
\def\theequation{5.3}\begin{equation}
\sum_{0\leq \ell \leq N}\, 
H_\ell \left(z,\bar t, 
\overline{b(t)} \right)\, 
a(t)^{N-\ell}=0,
\end{equation}
for all $t\in M \cap \Delta_n (\rho_3)$ and that the $\mathcal{
C}^\infty$-smooth functions $H_\ell \left(
z,\bar t, \overline{ b(t)} \right)$ are
not all identically zero. Then there exists an open neighborhood
$\mathcal{V}_4$ of $\mathcal{O}_{ CR}(M,p_0) \cap \Delta_n(\rho_3)$ in
$\C^n$ to which $a(t)$ extends holomorphically.
\end{mainproposition}

This proposition, applied to each $h_{ i'}$, completes
the proof of Theorem~2.9. The remainder of this section is
devoted to describe its proof, intuitively speaking.

\subsection*{5.4.~Heuristic} 
Although the equations~\thetag{3.10} (of a form similar
to~\thetag{ 5.3}) were obtained by applying the $(0,1)$ vector fields
tangent to $M$, it will be crucial to reapply the vector fields
$\overline{L}_k$ to~\thetag{ 5.3}.

But before applying the $\overline{ L}_k$, it is also crucial to
assume that $N$ is the smallest possible integer with the property
that there exists a relation of the form~\thetag{5.3} on $M\cap
\Delta_n(\rho_3)$. This assumption is of course free.  Since $a(t)$ is
CR, it may be considered as a constant by the derivations
$\overline{L}_k$.  Suppose for a while that dividing by $H_0 \left(
z,\bar t, \overline{b(t)} \right)$ is allowed in some sense. Then
rewrite~\thetag{5.3} as follows:
\def\theequation{5.5}\begin{equation}
a(t)^N + \sum_{1\leq \ell \leq N}\, 
{H_\ell \left(z,\bar t, 
\overline{b(t)}\right)\over
H_0\left(z,\bar t, \overline{b(t)}\right)}\, 
a(t)^{N-\ell}=0.
\end{equation}
Applying now the derivations $\overline{L}_k$, the term
$\overline{L}_k\left[ a(t)^N \right]$ vanishes (crucial fact), which
yields the identities
\def\theequation{5.6}\begin{equation}
\footnotesize
\sum_{1\leq \ell \leq N}\,
{H_0 \left(z,\bar t, \overline{b(t)}\right)
\, \overline{L}_k\left[H_\ell\left(z,\bar t, 
\overline{b(t)}\right)\right]-H_\ell\left(z,\bar t, 
\overline{b(t)}\right)\, \overline{L}_k\left[
H_0\left(z,\bar t, \overline{b(t)}\right)\right]
\over
\left[H_0\left(z,\bar t, 
\overline{b(t)}\right)\right]^2}\, 
a(t)^{N-\ell}=0,
\end{equation}
for $k= 1,\dots, m$. Now, after chasing the unnecessary denominator,
observe that since the coefficients of the $\overline{L}_k$ are
holomorphic in $\left( z,\bar t \right)$, there exist new holomorphic
functions $H_{1,\ell}$, $1\leq \ell \leq N$, such that
the identity~\thetag{ 5.6} may be rewritten as
\def\theequation{5.7}\begin{equation}
\sum_{1\leq \ell \leq N}\, 
H_{1,\ell} \left( z,\bar t, \overline{b_1(t)} \right)\, 
a(t)^{N-\ell}=0,
\end{equation}
where of course
\def\theequation{5.8}\begin{equation}
\overline{b_1(t)}:=
J_{\bar t}^1\overline{b(t)}=
\left(\partial_{\bar t_i} \overline{b_j(t)}\right)_{
1\leq i\leq n, \,  1\leq j\leq \nu}.
\end{equation}
Setting $\nu_1:= \nu (n+1)$, the relation~\thetag{ 5.7} is totally
similar to~\thetag{5.3}, if some freedom is allowed about the number
of functions $\overline{ b}_j$.  However, the degree $N-1$ of the
relation~\thetag{ 5.7} is strictly less than the degree $N$ of the
relation~\thetag{ 5.3}.  {\it Because the degree $N$ was chosen to be
the smallest possible one, this relation~\thetag{ 5.7} 
has to be trivial}.

Equivalently:
\def\theequation{5.9}\begin{equation}
\overline{L}_k\left({H_\ell \left( z,\bar t, 
\overline{b(t)} \right) \over
H_0 \left( z,\bar t, \overline{b(t)} \right)}
\right)\equiv 0,
\end{equation} 
over $M\cap \Delta_n (\rho_3)$, for $k=1, \dots,m$ and $\ell=1, \dots,
N$.  In other words, the quotient $H_\ell /H_0$ (which has to be
defined carefully in some sense) is CR on $M\cap \Delta_n(\rho_3)$.
Informally speaking, such an identity should be exceptional, because
the term $\overline{ b(t)}$ is anti-CR. In fact, one may expect
intuitively that if the relation~\thetag{5.9} is satisfied, then there
are no terms $\overline{ b(t)}$ at all, and hence the quotient
$H_\ell/ H_0$ extends meromorphically (but not holomorphically,
because of the presence of a quotient) to a neighborhood
$\mathcal{V}_3$ of $\mathcal{ O}_{ CR}(M, p_0) \cap \Delta_n (\rho_3)$
in $\C^n$.  This is true and will be proved below in Main Lemma~7.1
formulated below, where the assumption that $a(t)$ and $b(t)$ are
already real analytic in a neighborhood of the point $q_0 \in
\mathcal{ O}_{CR} (M,p_0)\cap \Delta_n(\rho_3)$ is used.

According to the works of K. Oka and E.E. Levi, a meromorphic function
defined in a domain of $\C^n$ is always the global quotient of two
holomorphic functions ({\it see}~\cite{ks} for instance). It follows
that there exist functions $R_\ell= R_\ell(t)$ holomorphic in
$\mathcal{ V}_3$ such that
\def\theequation{5.10}\begin{equation}
{H_\ell\left(z,\bar t, 
\overline{b(t)}\right)\over
H_0\left( z,\bar t, \overline{b(t)}\right)}\equiv
{R_\ell(t)\over R_0(t)},
\end{equation}
for all $t\in \mathcal{V}_3$. Replacing then~\thetag{ 5.10}
in~\thetag{5.5} and chasing the denominator:
\def\theequation{5.11}\begin{equation}
R_0(t)\, a(t)^N + \sum_{1\leq \ell \leq N}\, 
R_\ell(t)\, 
a(t)^{N-\ell}=0,
\end{equation}
for all $t\in \mathcal{ V}_3$. Here, the function $R_0 (t)$ is not
identically zero. Finally, by reproducing Section~3.3 of~\cite{ da} or
Proposition~6.4 of~\cite{ mmz1} ({\it see}
also~\cite{ n}), it follows that $a(t)$ is real
analytic at every point of $\mathcal{ O}_{CR} (M,p_0) \cap \Delta_n
(\rho_3)$, hence (thanks to the
Severi-Tomassini theorem)
it extends holomorphically to an open neighborhood
$\mathcal{ V}_4\subset \mathcal{V}_3$ of $\mathcal{ O}_{CR}( M,p_0)
\cap \Delta_n(\rho_3)$ in $\C^n$.

Sections~6, 7,8,9 and 10 are devoted to complete rigorously all the
details of this strategy of proof.

\section*{\S6.~Rings of CR functions and their fields of fractions}

Inspired by the preceding discussion, introduce the set $\mathcal{H}
\left( M,\rho_3, \varepsilon \right)$ of functions of the form
$H\left( z,\bar t, \overline{ b(t)} \right)$, where 
$\nu\in\N$, where $b(t)$ is a
$\C^\nu$-valued $\mathcal{ C}^\infty$-smooth CR mapping defined on $M
\cap \Delta_n (\rho_3)$ which extends holomorphically to $\mathcal{
W}_3$ and which is real analytic at the point $q_0 \in \mathcal{O}_{CR}
(M,p_0)\cap \Delta_n (\rho_3)$ and where $H= H \left(
z,\bar t, \bar b\right)$ is
holomorphic for $\vert z \vert < \rho_3$, $\vert \bar t \vert <
\rho_3$ and for $\left
\vert \bar b - \overline{ b(0)} \right\vert < \varepsilon$
and where $\left
\vert \overline{b(t)}- \overline{ b(0)} \right
\vert <
\varepsilon$ for all $t \in M \cap \Delta_n( \rho_3)$.  The main
feature of $\mathcal{ H} (M,\rho_3, \varepsilon)$ is the following.

\def\thelemma{6.1}\begin{lemma} 
Every function $H\left( z,\bar t, \overline{b(t)} \right)\in
\mathcal{H}(M,\rho_3,\varepsilon)$ admits a real analytic extension to
$\mathcal{W}_3$ which is antiholomorphic with respect to the complex
transversal coordinates $w=(w_1,\dots, w_d)\in\C^d$.
\end{lemma}

\proof 
Indeed, $\overline{b(t)}$ admits an antiholomorphic extension
to $\mathcal{W}_3$ by assumption, the function $H
\left( z,\bar t, \overline{b(t)} \right)$ 
is holomorphic with
respect to its variables, but it also depends on the holomorphic
variables $z=(z_1,\dots,z_m)$ in general. Consequently, the
antiholomorphicity with respect to $z$ (only) is lost.
\endproof

\def\thelemma{6.2}\begin{lemma}
The set $\mathcal{H}(M,\rho_3,\varepsilon)$ is
an {\rm entire} ring which is stable under differentiation 
by the $(0,1)$ vector fields tangent to $M$:
\def\theequation{6.3}\begin{equation}
\overline{L}_k \, \mathcal{H}(M,\rho_3,\varepsilon)\subset
\mathcal{H}(M,\rho_3,\varepsilon), \ \ \ \ \
k=1,\dots,m.
\end{equation}
\end{lemma}

\proof
The set $\mathcal{H}(M,\rho_3,\varepsilon)$ is obviously
stable under addition and multiplication. Suppose that
there exist two elements $H_1$, $H_2$ such that
\def\theequation{6.4}\begin{equation}
H_1 \left(z,\bar t, 
\overline{b_1(t)} \right)\cdot
H_2 \left(z,\bar t, 
\overline{b_2(t)} \right) = 0,
\end{equation}
for all $t\in M \cap \Delta_n(\rho_3)$. Then there exists a nonempty
open subset $V$ of $M \cap \Delta_n (\rho_3)$ on which $H_1( z, \bar
t, \overline{ b_1(t)}) $ or $H_2
\left( z, \bar t, \overline{ b_2(t)} \right)$\,--\,say
$H_1
\left( z, \bar t, \overline{ b_1(t)} \right)$ to fix 
ideas\,--\,vanishes identically. One must show that
that $H_1\left( z,\bar t, \overline{b(t)} \right)$ 
vanishes identically on
$M\cap \Delta_n(\rho_3)$. 

According to Lemma~6.1, the function $H_1 \left( 
z, \bar t, \overline{
b_1(t)} \right)$ 
extends real analytically and antiholomorphically with
respect to $w$ into  the wedge $\mathcal{W}_3$. 
By the principle of analytic
continuation (for real analytic functions), it suffices to show that $H_1 \left(
z, \bar t, \overline{ b_1(t)} \right)$ 
vanishes identically on a nonempty
open subset of $\mathcal{W}_3$.

Let $p\in V$ and let $\mathcal{V}_p$ be an open polydisc centered at
$p$ with $\mathcal{V}_p\cap M\subset V$.  For $q= (z_q,w_q)\in
\mathcal{V}_p$, the intersection $M\cap \{z=z_q\} \cap \mathcal{V}_p$
is maximally real in the slice $\{z=z_q\}\cap \mathcal{V}_p$. Also,
the function $H_1\left( z, \bar t, \overline{ b_1(t)} \right)$ extends
antiholomorphically with respect to $w$ into the sliced wedge
$\mathcal{W}_3 \cap \{z=z_q\} \cap \mathcal{V}_p$.  By the generic
uniqueness principle (for antiholomorphic functions), it follows that
$H_1 \left( z, \bar t, \overline{ b_1(t)}\right)$ vanishes identically
in the sliced wedge $\mathcal{W}_3 \cap \{z=z_q\} \cap
\mathcal{V}_p$. Since $z_q$ was arbitrary, it follows that $H_1
\left( z,
\bar t, \overline{ b_1(t)}\right)$
vanishes identically in the nonempty open
subset $\mathcal{W}_3 \cap \mathcal{V}_p$ of $\mathcal{W}_3$, as
desired.

Finally, since the coefficients of the vector fields $\overline{L}_k$
are holomorphic for $\vert z \vert < \rho_1$ and $\vert \bar t \vert <
\rho_1$, the second property follows from the chain rule.  The proof
of Lemma~6.2 is complete.
\endproof

A slightly more precise result has in fact been established.

\def\thecorollary{6.5}\begin{corollary}
The zero-set of a nonzero function $H
\left( z,\bar t,\overline{b(t)} \right)$ in
$\mathcal{H} \left( M,\rho_3,\varepsilon \right)$ is a
closed subset of $M\cap \Delta_n (\rho_3)$ with nonempty
interior. \qed
\end{corollary}

In the sequel, denote by $\mathcal{Z}_H$ the zero set of $H
\left( z,\bar t,
\overline{b(t)} \right)$ on $M\cap \Delta_n(\rho_3)$.  Since the ring
$\mathcal{H}(M,\rho_3,\varepsilon)$ is entire, it is allowed to
consider its {\sl field of fractions}
$\mathcal{R}(M,\rho_3,\varepsilon)$ which consists formally of
quotients of the form
\def\theequation{6.6}\begin{equation}
{H_1\left(z,\bar t, \overline{b_1(t)}\right)\over
H_2\left(z,\bar t, \overline{b_2(t)}\right)},
\end{equation}
and which may be viewed as a standard complex-valued function on the
dense open subset $(M\cap \Delta_n(\rho_3))\backslash
\mathcal{Z}_{H_2}$.

\subsection*{6.7.~Two notions of algebraic dependence}
Since $\mathcal{ R} (M, \rho_3, \varepsilon)$ is a field, the notion
of algebraic dependence makes sense. Precisely, a
$\mathcal{ C}^\infty$-smooth CR function $a(t)$ which is defined on
$M\cap \Delta_n (\rho_3)$ is called 
{\sl algebraic over} $\mathcal{ R}(M,
\rho_3, \varepsilon)$ if there exists a nonzero polynomial
\def\theequation{6.8}\begin{equation}
\sum_{0\leq \ell \leq N}\, 
{H_{1,\ell} \left(z,\bar t, 
\overline{b_{1,\ell}(t)} \right)\over
H_{2,\ell} \left(z,\bar t, 
\overline{b_{2,\ell}(t)} \right)}
\ X^{N-\ell},
\end{equation}
with coefficients in $\mathcal{R}(M, \rho_3,\varepsilon)$ which
annihilates $a(t)$, {\it i.e.} such that 
\def\theequation{6.9}\begin{equation}
\sum_{0\leq \ell \leq N}\, 
{H_{1,\ell}\left(z,\bar t, 
\overline{b_{1,\ell}(t)}\right)\over
H_{2,\ell}\left(z,\bar t, 
\overline{b_{2,\ell}(t)}\right)}
\ a(t)^{N-\ell}=0,
\end{equation}
for all $t\in (M\cap \Delta_n(\rho_3)) \backslash \bigcup_{0\leq \ell
\leq N}\, \mathcal{Z}_{H_{2, \ell}}$.

Equivalently, after chasing the denominators, the
$\mathcal{C}^\infty$-smooth CR function $a(t)$ is algebraically
dependent over the ring $\mathcal{H}(M,\rho_3,\varepsilon)$.

Another notion of algebraic dependence is the following. Let
${\rm Mer}(\mathcal{O}_{CR}(M,p_0),
\rho_3)$ denote the field
of functions which are meromorphic in
some connected open neighborhood of $\mathcal{O}_{CR}(M,p_0)
\cap \Delta_n(\rho_3)$ in $\C^n$. In
this field, the classical algebraic operations are defined up to some
shrinking of the domains of definition. Then a
$\mathcal{C}^\infty$-smooth CR function $a(t)$ defined on $M\cap
\Delta_n(\rho_3)$ is called 
{\sl algebraic over} ${\rm Mer}(\mathcal{O}_{CR}(M,p_0),
\rho_3)$ if
there exists a nonzero polynomial
\def\theequation{6.10}\begin{equation}
\sum_{0\leq \ell \leq N}\, 
R_\ell (t) \
X^{N-\ell},
\end{equation}
with coefficients $R_\ell(t)\in{\rm Mer} (\mathcal{O}_{CR}(M,p_0),
\rho_3)$ which annihilates
$a(t)$, {\it i.e.} such that 
\def\theequation{6.11}\begin{equation}
\sum_{0\leq \ell \leq N}\, 
R_\ell (t) \
a(t)^{N-\ell}=0,
\end{equation}
for all $t$ in some neighborhood 
$\mathcal{V}_3$ of $\mathcal{O}_{CR}(M,p_0)\cap \Delta_n(\rho_3)$ in 
$\C^n$, outside the union of the polar sets
of the $R_\ell(t)$.

It is now time to reformulate and to slightly generalize the heuristic
discussion of \S5.4.  In the following lemma, a Main Lemma~7.1 (to be
formulated and to be proved later), is hidden.

\def\thelemma{6.12}\begin{lemma}
Let $a(t)$ be a $\mathcal{C}^\infty$-smooth CR function on $M\cap
\Delta_n(\rho_3)$ which extends holomorphically to the complete wedge
$\mathcal{W}_3$ in $\Delta_n(\rho_3)$
with edge $M\cap \Delta_n(\rho_3)$. Assume that $a(t)$ is
real analytic at one point $q_0\in \mathcal{O}_{CR}(M,p_0)\cap 
\Delta_n(\rho_3)$. If the
function $a(t)$ is algebraic over the field
$\mathcal{R}(M,\rho_3,\varepsilon)$ {\rm(}or 
equivalently over the ring $\mathcal{H}(M,\rho_3,\varepsilon)${\rm )}, 
then it is algebraic over the
field ${\rm Mer}(\mathcal{O}_{CR}(M,p_0),\rho_3)$.
\end{lemma}

\proof
Without loss of generality, assume that $a(t)$
satisfies a polynomial relation
\def\theequation{6.13}\begin{equation}
\sum_{0\leq \ell \leq N}\, 
H_\ell \left( z,\bar t, \overline{b(t)} \right)\
a(t)^{N-\ell}=0,
\end{equation}
for all $t\in M\cap \Delta_n(\rho_3)$,
{\it whose degree $N\geq 1$ is the
smallest possible}. Hence $H_0
\left( z,\bar t, \overline{b(t)} \right)$ does not
vanish identically. Proceeding exactly 
as in the heuristic discussion of
\S5.4, divide by $H_0
\left( z,\bar t, \overline{b(t)} \right)$, which yields
\def\theequation{6.14}\begin{equation}
a(t)^N + \sum_{1\leq \ell \leq N}\, 
{H_\ell \left( z,\bar t, 
\overline{b(t)} \right) \over
H_0 \left( z,\bar t, \overline{b(t)} \right)}\, 
a(t)^{N-\ell}=0.
\end{equation}
Apply now the derivations $\overline{L}_k$, which yields
\def\theequation{6.15}\begin{equation}
\footnotesize 
\sum_{1\leq \ell \leq N}\, 
{H_0\left( z,\bar t, \overline{b(t)} \right)
\, \overline{L}_k\left[H_\ell \left(z,\bar t, 
\overline{b(t)} \right)\right]-H_\ell \left(z,\bar t, 
\overline{b(t)} \right)\, \overline{L}_k\left[
H_0 \left(z,\bar t, \overline{b(t)} \right)\right]
\over
\left[H_0 \left(z,\bar t, \overline{b(t)} 
\right) \right]^2}\, 
a(t)^{N-\ell}=0.
\end{equation}
Now, after chasing the unnecessary denominator, and
using~\thetag{6.3}, rewrite this identity 
under the form
\def\theequation{6.16}\begin{equation}
\sum_{1\leq \ell \leq N}\, 
H_{1,\ell}\left( z,\bar t, \overline{b_1(t)} \right)\, 
a(t)^{N-\ell}=0,
\end{equation}
where $H_{1,\ell} \left(
z,\bar t, \overline{b_1(t)} \right)$ belongs to
$\mathcal{H}(M,\rho_3,\varepsilon)$. Thus a
relation which is totally similar to~\thetag{10.1.47}, but which is of
degree strictly less than $N$, has been
constructed. Because the degree $N$ was chosen to be
the smallest possible one, this relation has to be trivial.
Consequently, the $(0,1)$ vector fields $\overline{ L}_k$ annihilate
the quotients $H_\ell /H_0$ outside the zero-set
$\mathcal{Z}_{H_0}$. By Main Lemma~7.1 
(to be discussed later), this implies that for each $\ell=1,\dots,N$,
the quotient $H_\ell/H_0$ coincides in some neighborhood
$\mathcal{V}_3$ of $\mathcal{O}_{CR}(M,p_0)\cap \Delta_n(\rho_3)$ in
$\C^n$ with the restriction to $[\mathcal{O}_{CR}(M,p_0)\cap
\Delta_n(\rho_3)] \backslash \mathcal{Z}_{H_0}$ of a meromorphic
function $R_\ell(t)/ R_0(t)$. 
Finally, replacing $H_\ell/H_0$ by $R_\ell
/ R_0(t)$ in
the original relation~\thetag{6.13}, it follows that $a(t)$ is
algebraic over ${\rm Mer}(\mathcal{O}_{CR}(M,p_0),\rho_3)$.
The proof of Lemma~6.12 is complete.
\endproof

Thanks to this Lemma~6.12, the proof of Main 
Proposition~5.2 is 
complete, as explained after~\thetag{ 5.11}.
\qed

\section*{\S7.~Statement of Main Lemma~7.1}

In sum, everything relies upon a statement,
interesting in itself, which is
reformulated carefully, including
all the assumptions.

\def\themainlemma{7.1}\begin{mainlemma}
As in Main Proposition~5.2, let $M$ be a real analytic local
generic submanifold defined
in $\Delta_n(\rho_1)$, as in Definition~2.2
{\bf (I)}. Let $\rho_3$ with $0< \rho_3
< \rho_2< \rho_1$ and let $\mathcal{W}_3$ be a 
complete wedge in $\Delta_n( \rho_3)$ 
with edge $M\cap
\Delta_n(\rho_3)$. Let $\nu_1, \nu_2\in\N$, let $\varepsilon>0$ and
let $b_1(t), b_2(t)$ be two $\C^{\nu_1}$- and $\C^{\nu_2}$-valued
$\mathcal{C}^\infty$-smooth CR mappings defined on $M\cap
\Delta_n(\rho_3)$ which satisfy $\vert b_1(t)- b_1(0) \vert <
\varepsilon$, 
$\vert b_2(t)- b_2(0) \vert < \varepsilon$, which extend
both holomorphically to $\mathcal{W}_3$ and which are real analytic in
a neighborhood of the same point $q_0\in \mathcal{O}_{CR}(M,p_0) \cap
\Delta_n (\rho_3)$.  Let $H_1(z,\bar t, \overline{b}_1)$ and
$H_2\left(z,\bar t, \overline{b}_2\right)$ be two holomorphic
functions defined for $\vert z\vert < 
\rho_3$, for $\vert \bar t \vert
< \rho_3$ and for $\left \vert \overline{b}_1 -
\overline{b_1(0)}\right\vert < \varepsilon$, $\left\vert
\overline{b}_2 - \overline{b_2(0)}\right\vert < \varepsilon$. Assume
that $H_2\left( z,\bar t, \overline{b_2(t)} \right)$ does not vanish
identically on
$M\cap \Delta_n(\rho_3)$ and consider the quotient
\def\theequation{7.2}\begin{equation}
{H_1 \left(z,\bar t, \overline{b_1(t)}\right)\over
H_2\left(z,\bar t, \overline{b_2(t)}\right)},
\end{equation}
which belongs to $\mathcal{R}(M,\rho_3,\varepsilon)$ and remind that
by Corollary~6.5, the zero-set $\mathcal{Z}_{H_2}$ of $H_2
\left(z,\bar
t, \overline{b_2(t)}
\right)$ is a closed subset of $M\cap \Delta_n(\rho_3)$
with nonempty interior.  Then the following three conditions are
equivalent:
\begin{itemize}
\item[{\bf (1)}]
there exists a meromorphic function in 
${\rm Mer}(\mathcal{O}_{CR}(M,p_0),\rho_3)$ of
the form $P_1(t)/P_2(t)$, where $P_1(t)$ and $P_2(t)$ are holomorphic
functions defined in some open connected neighborhood $\mathcal{V}_3$
of $\mathcal{O}_{CR}(M,p_0)
\cap \Delta_n(\rho_3)$ in $\C^n$ with $P_2(t)\not \equiv 0$ such that
\def\theequation{7.3}\begin{equation}
{H_1\left(z,\bar t, \overline{b_1(t)}\right)\over
H_2\left(z,\bar t, \overline{b_2(t)}\right)}=
{P_1(t)\over P_2(t)},
\end{equation}
for all $t\in \mathcal{V}_3\cap M$ outside the zero
set of $H_2$;
\item[{\bf (2)}]
the quotient $H_1/H_2$ is CR over the dense open subset 
$(M\cap \Delta_n(\rho_3)) \backslash \mathcal{Z}_{H_2}$;
\item[{\bf (3)}]
There exists a nonempty open subset $V$ of the dense open subset
$(M\cap \Delta_n(\rho_3)) \backslash \mathcal{Z}_{H_2}$ on which the
quotient $H_1/H_2$ is CR.
\end{itemize}
\end{mainlemma}

\proof
Obviously, {\bf (2)} $\Rightarrow$ {\bf (3)}.  Treat first the two
implications {\bf (1)} $\Rightarrow$ {\bf (2)} and {\bf(3)}
$\Rightarrow$ {\bf (2)}, which are easy. Indeed, by Lemma~6.1, 
for $k=1, \dots, m$, the CR derivatived functions 
$\overline{ L}_k ( H_1/
H_2)$ also belong to $\mathcal{ R}(M, \rho_3, \varepsilon)$. If {\bf
(1)} or {\bf (3)} holds, then $H_1/ H_2$ is CR on a nonempty open
subset $V$ of $(M\cap \Delta_n(\rho_3)) \backslash \mathcal{Z}_{H_2}$,
so Corollary~6.5 implies that $\overline{ L}_k (H_1/H_2)$ vanishes
identically on $M\cap \Delta_n(\rho_3)\backslash \mathcal{Z}_{H_2}$,
which yields {\bf (2)}.

Next, begin the proof of the implication {\bf (2)} $\Rightarrow$
{\bf (1)}, which is by far the main task. Even if it is already 
known that {\bf (2)} and {\bf (3)} are equivalent, 
property {\bf (1)} will be established assuming only that 
$H_1/H_2$ is CR over a nonempty open subset $V$, as
in {\bf (3)}.

Define the Schwarz reflection across $M$ which 
stabilizes the slices $\{z=ct.\}$ explicitely by
\def\theequation{7.4}\begin{equation}
s: (z,w)\mapsto \left(z,\overline{\Theta}(z,\bar z,\bar w) \right).
\end{equation}
Shrinking $\rho_3$ a bit if necessary, one can construct a slightly
smaller complete wedge $\widetilde{ \mathcal{ W}}_3$ which is contained
in $s(\mathcal{W}_3)$, as depicted in {\sc Figure~1} above.
What is important is that $\mathcal{W}_3$ and 
$\widetilde{\mathcal{W}}_3$ are directed by two opposite
vectors at every point of $M\cap \Delta_n(\rho_3)$, because
a version of the edge of the wedge theorem 
will be applied in the end of the proof.

Since $s$ is antiholomorphic with respect to $w$ and since 
$H_1\left(z,\bar
t, \overline{b_1 (t)}\right)$ 
and $H_2\left(z, \bar t, \overline{ b_2(t)}\right)$ extend
as real analytic functions in $\mathcal{W}_3$ which are
antiholomorphic with respect to $w$, it follows by composition with
the reflection $s(z,w)$ that $H_1$ and $H_2$ extend to be real
analytic in $\widetilde{\mathcal{ W}}_3$ and holomorphic with respect
to $w$. Denote by $\widetilde{ H}_1 (t,\bar t)$ and $\widetilde{ H}_2
(t,\bar t)$ these two real analytic extensions. Since the
antiholomorphic reflection coincides with the identity mapping on $M
\cap \Delta_n( \rho_3 )$, the $\mathcal{ C }^\infty$ boundary values
of the real analytic extensions $\widetilde{ H}_1$, $\widetilde{ H}_2$
are just $H_1$ and $H_2$.

Since $\widetilde{H}_1/ \widetilde{H }_2$ is CR on the nonempty open
subset $V$ of $(M \cap \Delta_n(\rho_3 ))\backslash \mathcal{Z}_2$,
and since it is holomorphic with respect to the variable $w$ in the
vertical conelike slices $\{z=ct.\}\cap \widetilde{ \mathcal{W}}_3$,
it follows from an elementary (and known)
separate Cauchy-Riemann principle
that there exists a neighborhood $\mathcal{ V}$ of $V$ in $\C^n$ such
that $\widetilde{ H}_1/ \widetilde{H}_2$ extends holomorphically to
$\mathcal{ V}\cap \widetilde{ \mathcal{W}}_3$. Apply now the
following general lemma to deduce that $\widetilde{ H}_1/
\widetilde{H}_2$ can be represented as a quotient $P_1/P_2$, where
$P_1$ and $P_2\neq 0$ are holomorphic in $\widetilde{\mathcal{W}}_3$.

\def\thelemma{7.5}\begin{lemma}
Let $\Omega \subset 
\C^n$ be a domain, let $A_1$ and $A_2\neq 0$ be two
real analytic functions in $\Omega$. Suppose that there exist a point
$p\in \Omega\backslash \mathcal{Z}_{A_2}$ and a nonempty open
neighborhood of $p$ such that $V_p$ is contained in $\Omega
\backslash \mathcal{Z}_{A_2}$ and such that the restriction to $V_p$ of
$A_1/A_2$ is holomorphic in $V_p$. Then there exist two holomorphic
functions $P_1$ and $P_2\neq 0$ in $\Omega$ such that 
\def\theequation{7.6}\begin{equation}
\left.{P_1\over P_2}\right\vert_{\Omega 
\backslash \mathcal{Z}_{A_2}}
=
\left.{A_1\over A_2}\right\vert_{\Omega 
\backslash \mathcal{Z}_{A_2}}.
\end{equation}
\end{lemma}

\proof
In this proof (only), denote some complex coordinates on $\C^n$ by
$z=(z_1, \dots, z_n)$ and $z= x+iy$.
By assumption, for $i=1,\dots,n$ and for
$z=x+iy\in V_p$, it holds that
\def\theequation{7.7}\begin{equation}
{\partial \over \partial \bar z_i} \left(
{A_1(x,y)\over A_2(x,y)}\right)\equiv 0,
\end{equation}
which yields by analytic continuation:
\def\theequation{7.8}\begin{equation}
A_2(x,y) \, {\partial A_1\over \partial \bar z_i}(x,y)-
A_1(x,y) \, {\partial A_2\over \partial \bar z_i}(x,y)
\equiv 0,
\end{equation}
for all $z\in \Omega$. It follows that the
quotient $A_1/A_2$ is holomorphic on $\Omega\backslash
\mathcal{Z}_{A_2}$, even if $\Omega\backslash
\mathcal{Z}_{A_2}$ is not (locally) connected.
In fact, to achieve the proof, the local connectedness of
$\Omega\backslash \mathcal{Z}_{A_2}$ is needed.
To fix this point, it is sufficient to show that $\mathcal{Z}_{A_2}$
is in fact of real dimension $\leq 2n-2$, hence in particular
$\Omega\backslash \mathcal{Z}_{A_2}$ is locally connected.

Indeed, proceeding by contradiction, in a neighborhood of a point
$q\in \mathcal{Z}_{A_2}$ at which $\mathcal{Z}_{A_2}$ is geometrically
smooth and of dimension $2n-1$, make a translation of
coordinates so that $q$ is the origin and $\mathcal{Z}_{A_2}$ is
locally represented by $\{x_1=0\}$. Then write
\def\theequation{7.9}\begin{equation}
A_2(x,y)= (x_1)^\kappa [A_{2,0}(
x_2, \cdots,x_n,y_1,\dots,y_n)+x_1\, A_{2,1}(x,y)],
\end{equation}
for some integer $\kappa \geq 1$ and some real analytic function
$A_{2,0}$ independent of $x_1$ which is not identically zero. Without
loss of generality, assume that $A_1(x,y)\not \equiv 0$, hence
\def\theequation{7.10}\begin{equation}
A_1(x,y)= (x_1)^\lambda [A_{1,0}(
x_2, \cdots,x_n,y_1,\dots,y_n)+x_1\, A_{1,1}(x,y)],
\end{equation}
for some integer $\lambda \geq 0$, where the function $A_{1,0}$ is
also not identically zero. Since by assumption $q \in \mathcal{ Z}_{
A_2}$, it holds that $\kappa \geq \lambda+1$.
Compute then~\thetag{7.8} for $i = 1$, which yields
\def\theequation{7.11}\begin{equation}
0\equiv 
\left( \frac{ \lambda - \kappa} {2}\right)  \, 
x_1^{\kappa + \lambda-1}[A_{1,0} \, A_{2,0} + {\rm O}(x_1)],
\end{equation}
a contradiction.

Thus $A_1/A_2$ extends holomorphically to the {\it
locally connected}\, open set $\Omega\backslash 
\mathcal{ Z}_{A_2}$.

To pursue the proof, choose an arbitrary point $q\in \Omega$, 
center the coordinates $z$ at $q$
and consider the polydisc $\Delta_n (\rho_q ) \subset \Omega$,
where $\rho_q= \inf_{ r\in \partial \Omega }\,\vert r-q \vert$.  Then
$A_1(x,y)$ and $A_2(x,y)$ may be developed in power series with
respect to $x$ and $y$ with radius of convergence at least equal to
$\rho_q$. After perhaps making a unitary transformation, 
assume that the maximally real submanifold $[\R^n\times \{0\}]\cap
\Delta_n(\rho_q)$ is not contained in $\mathcal{Z}_{A_2}$. Then
complexify the variable $x\in \R^n$ to be the complex variable $z\in
\C^n$ and introduce the element of meromorphic function
\def\theequation{7.12}\begin{equation}
R_q(z):=
{A_1(z,0)\over A_2(z,0)},
\end{equation}
which is
defined in $\Delta_n(\rho_q)$. By construction,
the restrictions to the maximally real subspace 
\def\theequation{7.13}\begin{equation}
[(\R^n\times \{0\})\cap 
\Delta_n (\rho_q)] \backslash \mathcal{Z}_{A_2}
\end{equation}
of $R_q (z)$ and of $A_1(x,y) /A_2 (x,y)$ coincide. Because $\Delta_n
(\rho_q) \backslash \mathcal{Z}_{A_2}$ is locally connected, it
follows from the principle of analytic continuation that $R_q (z)$ and
$A_1 (x,y) / A_2 (x,y)$ coincide all over $\Delta_n( \rho_q)
\backslash \mathcal{ Z}_{A_2}$.

In summary, at every point $q$ of $\Omega$, a
meromorphic extension of the holomorphic function $A_1(x,y)/A_2(x,y)$
defined in $\Omega \backslash \mathcal{Z}_{A_2}$ has
been constructed. It follows
from~\cite{ ks} that this meromorphic function can be represented as a
quotient of holomorphic functions in 
$\Omega$, which completes the proof of the
lemma.
\endproof

\subsection*{7.14.~Continuation of the proof}
Thus $\widetilde{ H_1} / \widetilde{ H_2}$ 
extends meromorphically to $\widetilde{ \mathcal{ W}}_3
\cup \mathcal{ V}$. Of course, neither $\mathcal{ V}$
nor a neighborhood (in $\C^n$) of the
dense subset $(M \cap \Delta_n (\rho_3)) 
\backslash \mathcal{ Z}_{ H_2}$ need contain a
point of the (thin in the nonminimal case) subset
$\mathcal{ O}_{ CR} (M, p_0) \cap 
\Delta_n (\rho_3)$. 

Here comes the assumption 
that the CR mappings $b_1( t)$ and $ b_2(t)$
extend holomorphically to a neighborhood 
$\mathcal{ V}_{ q_0}$ (in $\C^n$) of some point $q_0 \in 
\mathcal{ O}_{ CR} (M, p_0) \cap 
\Delta_n (\rho_3)$. It follows that
$H_1 \left( z, \bar t, \overline{ b_1( t)} \right)$ and 
$H_2 \left( z, \bar t, \overline{ b_2( t)} \right)$
extend real analytically to $\mathcal{ V}_{q_0}$.
Thanks to the (already established)
equivalence between {\bf (2)} and {\bf (3)}, 
the quotient $H_1 / H_2$ is CR over the dense
open subset $(M \cap \Delta_n (\rho_3) )
\backslash \mathcal{ Z}_{ H_2}$.
In particular, there exists a point
$r_0 \in \mathcal{ V}_{q_0}$ in a neighborhood of
which $H_2$ is nonzero. Then, thanks to the
Severi-Tomassini extension theorem, $H_1 / H_2$ 
extends holomorphically to a neighborhood
$\mathcal{ V}_{ r_0}$ (in $\C^n$) of $r_0$.
Applying Lemma~7.5 just above, it follows that 
$H_1/H_2$ extends meromorphically to 
$\mathcal{ V}_{ q_0}$. Forget the
open subset $\mathcal{ V}$ and summarize the
obtained extension result.

\def\thelemma{7.15}\begin{lemma}
The CR quotient $H_1/H_2$ extends meromorphically to 
$\widetilde{\mathcal{W}}_3\cup \mathcal{V}_{q_0}$.
\end{lemma}

The goal, to which the remainder of the paper is devoted, is to prove
that $H_1/H_2$ extends meromorphically to some open neighborhood
$\mathcal{V}_3$ of $\mathcal{O}_{CR}(M,p_0)\cap \Delta_n(\rho_3)$ in
$\C^n$.

To this aim, define the set $D$ of points $q\in \mathcal{O}_{CR}
(M,p_0)\cap \Delta_n( \rho_3)$ such that there exists a small nonempty
open polydisc $\mathcal{ V}_q$ centered at $q$ and a meromorphic
extension $R_q (t)$ of $H_1/H_2 \vert_{ ( \mathcal{V}_q \cap M)
\backslash \mathcal{ Z}_{ H_2}}$ to $\mathcal{ V}_q$. This set is
nonempty, since $q_0$ belongs to $D$ by assumption. It follows from
the uniqueness principle on a generic edge that the various
meromorphic functions $R_q(t)$ glue together to provide a well defined
meromorphic function $R_D(t)$ defined in the open neighborhood
$\mathcal{V}_D:= \bigcup_{ p\in D}\, \mathcal{V}_q$ of $D$ in $\C^n$.
State this property as a step lemma.

\def\thelemma{7.16}\begin{lemma}
The CR quotient extends meromorphically to 
$\widetilde{\mathcal {W}}_3 \cup \mathcal{V}_D$.
\end{lemma}

If $D = \mathcal{O}_{ CR}( M,p_0) \cap \Delta_n (\rho_3)$, Main
Lemma~7.1 would be proved, almost gratuitously. Suppose therefore that
the complement $E$ of $D$ in $\mathcal{O}_{CR}(M,p_0)\cap
\Delta_n(\rho_3)$ is nonempty. To conclude the proof of Main
Lemma~7.1, it will suffice to derive a contradiction in the following
form: establish that there exists in fact {\it at least one point}
$p_1\in E$ at which $H_1/ H_2$ extends meromorphically.

\section*{\S8.~Localization at a nice boundary point}

So, assume that $E$ (the bad set) and $D$ (the good set) 
are nonempty,
with $E\cap D= \emptyset$ and with $E\cup D=\mathcal{O}_{
CR}(M,p_0)\cap \Delta_n( \rho_3)$.  For technical convenience, it is
better to pick a special point $p_1 \in E$ so that $E$ lies behind a
real analytic generic ``wall'' $M_1$ passing through $p_1$, as
depicted in {\sc Figure~2} below
({\it cf.}~\cite{ mp1}, \cite{ mp2}).

\def\thelemma{8.1}\begin{lemma}
There exists a point $p_1\in E$ and a real analytic hypersurface
$M_1\subset M$ passing through $p_1$ which is {\rm generic} in $\C^n$
such that $E\backslash \{p_1\}$ lies, near $p_1$, in one side of
$M_1$.
\end{lemma}

\proof
Let $q\in E\neq\emptyset$ be an arbitrary point and let $\gamma$ be a
piecewise real analytic curve running in complex tangential directions
to $M$ (CR-curve) which links the point $q$ with the point $q_0\in D$.
Such a curve exists because the CR orbit $\mathcal{O}_{CR}(M,p_0)\cap
\Delta_n(\rho_3)$ is locally minimal at every point.  After shortening
$\gamma$ and changing the points $q$ and $q_0$ if necessary, one can
assume that $\gamma$ is a smoothly embedded segment, that $q$ and
$q_0$ belong to $\gamma$ and are close to each other. Therefore
$\gamma$ can be described as a part of an integral curve of some
nonvanishing real analytic section $L$ of $T^cM$ defined in a
neighborhood of $q$.

\bigskip
\begin{center}
\begin{picture}(0,0)%
\includegraphics{blowing-ellipsoids.pstex}%
\end{picture}%
\setlength{\unitlength}{3947sp}%
\begingroup\makeatletter\ifx\SetFigFont\undefined%
\gdef\SetFigFont#1#2#3#4#5{%
  \reset@font\fontsize{#1}{#2pt}%
  \fontfamily{#3}\fontseries{#4}\fontshape{#5}%
  \selectfont}%
\fi\endgroup%
\begin{picture}(5799,2432)(1211,-2150)
\put(1947, 81){\makebox(0,0)[lb]{\smash{{\SetFigFont{7}{8.4}{\familydefault}{\mddefault}{\updefault}{\color[rgb]{0,0,0}$L$}%
}}}}
\put(1306,-759){\makebox(0,0)[lb]{\smash{{\SetFigFont{7}{8.4}{\familydefault}{\mddefault}{\updefault}{\color[rgb]{0,0,0}$\gamma$}%
}}}}
\put(3171,-381){\makebox(0,0)[lb]{\smash{{\SetFigFont{7}{8.4}{\familydefault}{\mddefault}{\updefault}{\color[rgb]{0,0,0}$\Upsilon$}%
}}}}
\put(4116,-471){\makebox(0,0)[lb]{\smash{{\SetFigFont{7}{8.4}{\familydefault}{\mddefault}{\updefault}{\color[rgb]{0,0,0}$Q_{\delta_1}$}%
}}}}
\put(4946,-1004){\makebox(0,0)[lb]{\smash{{\SetFigFont{7}{8.4}{\familydefault}{\mddefault}{\updefault}{\color[rgb]{0,0,0}$p_1$}%
}}}}
\put(4671,-1339){\makebox(0,0)[lb]{\smash{{\SetFigFont{7}{8.4}{\familydefault}{\mddefault}{\updefault}{\color[rgb]{0,0,0}$M_1$}%
}}}}
\put(5566,-349){\makebox(0,0)[lb]{\smash{{\SetFigFont{7}{8.4}{\familydefault}{\mddefault}{\updefault}{\color[rgb]{0,0,0}$E$}%
}}}}
\put(6721,-789){\makebox(0,0)[lb]{\smash{{\SetFigFont{7}{8.4}{\familydefault}{\mddefault}{\updefault}{\color[rgb]{0,0,0}$\gamma$}%
}}}}
\put(3161,-900){\makebox(0,0)[lb]{\smash{{\SetFigFont{7}{8.4}{\familydefault}{\mddefault}{\updefault}{\color[rgb]{0,0,0}$q$}%
}}}}
\put(6156,-834){\makebox(0,0)[lb]{\smash{{\SetFigFont{7}{8.4}{\familydefault}{\mddefault}{\updefault}{\color[rgb]{0,0,0}$q_0$}%
}}}}
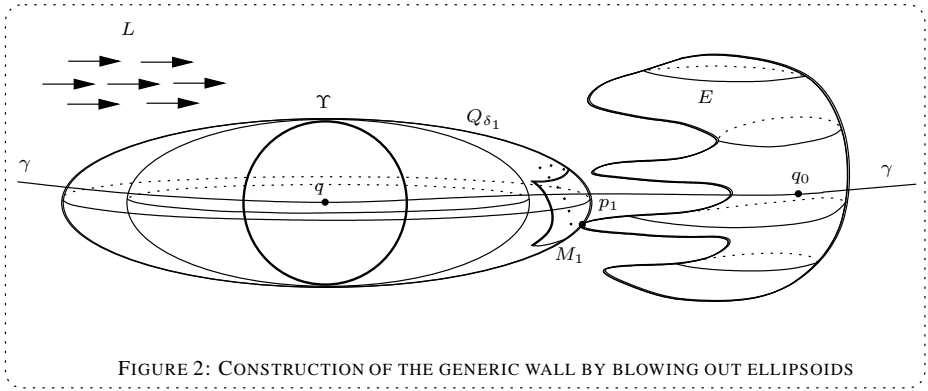
\put(1927,-2058){\makebox(0,0)[lb]{\smash{{\SetFigFont{8}{9.6}{\familydefault}{\mddefault}{\updefault}{\color[rgb]{0,0,0}{\sc Figure~2: Construction of the generic wall by blowing out ellipsoids}}%
}}}}
\end{picture}%

\end{center}
\bigskip

Let $H\subset M$ be a small $(2m+d-1)$-dimensional real analytic
hypersurface passing through $q$ such that $L(q)$ is not tangent to
$H$ at $q$. Integrating $L$ with initial conditions in $H$, 
one obtains
real analytic coordinates $(s_1,s_2)\in \R\times \R^{2m+d-1}$ so that
for fixed $s_{2,0}$, the segments $(s_1,s_{2,0})$ are contained in the
trajectories of $L$. After a translation, one may assume that the
origin $(0,0)$ corresponds to a point of $\gamma$ close to $q$ which
is not contained in $E$, again denoted by $q$. Fix a small
$\varepsilon>0$ and for real $\delta\geq 1$, define the ellipsoids
({\it see} {\sc Figure~2} above)
\def\theequation{8.2}\begin{equation}
Q_{\delta}:=\{(s_1,s_2): \, 
\vert s_1 \vert^2/\delta+\vert s_2 \vert^2< \varepsilon\}.
\end{equation}
There exists the 
smallest $\delta_1>1$ with $\overline{Q_{\delta_1}}\cap
E\neq \emptyset$.  Then $\overline{Q_{\delta_1}} \cap E= \partial
Q_{\delta_1} \cap E$ and $Q_{\delta_1} \cap E=\emptyset$.  Observe
that every boundary $\partial Q_{\delta}$ is transverse to the
trajectories of $L$ out off the equatorial set $\Upsilon:=\{(0,s_2):
\, \vert s_2 \vert^2=\varepsilon\}$ which is contained in $D$. Hence
$\partial Q_{\delta_1}$ is transverse to $L$ at every point of
$\partial Q_{\delta_1}\cap E$. So $\partial Q_{\delta_1} \backslash
\Upsilon$ is generic in $\C^n$, since $L$ is a section of $T^cM$.

To conclude, it suffices to choose a point $p_1\in \partial
Q_{\delta_1} \cap E$ and to take for $M_1$ a small real analytic
hypersurface passing through $p_1$ which is tangent to $\partial
Q_{\delta_1}$ and satisfies $M_1\backslash \{p_1\}\subset
Q_{\delta_1}$. The proof of Lemma~8.1. is complete.
\endproof

\subsection*{8.3.~Localization}
Choose now a point $p_1\in E$ as in Lemma~10.1 and
choose $\varepsilon_1>0$ such that the polydisc
$\Delta_n(\varepsilon_1)$ (in the new coordinates centered at $p_1$)
with center $p_1$ is contained in the polydisc $\Delta_n(\rho_3)$ (in
the old coordinates centered at $p_0$) with center $p_0$ and
localize everything in $\Delta_n(\varepsilon_1)$ ({\it cf.}~{\sc
Figure~3} below, where $E$ has been redrawn on the
left).

Denote by
$M_1^-$ the negative, left open one-sided neighborhood of $M_1$ in $M$
such that $E\backslash \{p_1\}$ lies in $M_1^-$ in a neighborhood of
$p_1$, and denote by $M_1^+$ the other side, which is contained in
$D$ by assumption. Choose affine coordinates vanishing at $p_1$,
still denoted by $t=(z,w)=(x+iy,u+iv)\in \C^m\times \C^d$ in
order that $T_0M=\{v=0\}$ and $T_0 M_1=\{v=0, x_1=0\}$. Denote
by $z_\sharp\in \C^{m-1}$ the coordinates $(z_2, \dots, z_m)$, 
which are
decomposed in real and imaginary part as $z_\sharp =x_\sharp
+iy_\sharp$.  Then $M$ is defined by $v_j =\varphi_j(x,y,u)$,
$j=1,\dots,d$ with $\varphi_j(0) =0$, $d\varphi_j(0)=0$ and 
$M_1$
is defined by a supplementary equation $x_1=\psi(y_1,x_\sharp,
y_\sharp, u)$, with $\psi(0)=0$, $d\psi(0)=0$.

After possibly replacing $M_1$ by a new hypersurface which is
contained in $M_1^+\cup \{p_1\}$ ({\it cf.} the end of the proof of
Lemma~8.1) and after possibly making a dilatation of coordinates, 
one may
assume that the supplementary equation of $M_1$ is simply given by
$x_1= y_1^2+\vert z_\sharp \vert^2+ \vert u \vert^2$, 
and that $M_1^+$ is given by
\def\theequation{8.4}\begin{equation}
M_1^+: \ \ \ \ \ 
x_1> y_1^2+ \vert z_\sharp \vert^2+ \vert u \vert^2.
\end{equation}
(For the readability of {\sc Figure~3} below,  
the curvature of $M_1$ has been reversed, hence the picture of
$M_1$ is slightly incorrect).

Intersect everything with the polydisc 
$\Delta_n(\varepsilon_1)$ centered at $p_1$
In particular, define 
\def\theequation{8.5}\begin{equation}
\widetilde{\mathcal{W}}_{
\varepsilon_1}:=
\widetilde{\mathcal{W}}_3 
\cap \Delta_n(\varepsilon_1) \ \ \ \ \ \
{\rm and} \ \ \ \ \ 
\mathcal{V}_{\varepsilon_1}:= 
\mathcal{V}_D\cap \Delta_n(\varepsilon_1).
\end{equation}

To reach the desired contradiction 
({\it cf.}~the last sentence
of \S7) which achieves the
proof of Main Lemma~7.1, 
the principal goal is now to show that the
quotient
\def\theequation{8.6}\begin{equation}
{H_1\left(z,\bar t, \overline{b_1(t)}\right)\over
H_2\left(z,\bar t, \overline{b_2(t)}\right)},
\end{equation}
which already extends holomorphically to
$\widetilde{W}_{\varepsilon_1} \cup \mathcal{V}_{\varepsilon_1}$
extends holomorphically to a neighborhood of $p_1$ in $\C^n$.  For
this purpose, the technique of normal deformations of analytic discs
enters the scene (\cite{ tu2}, \cite{ mp1}, 
\cite{ mp2}).

\section*{ \S9.~Construction of analytic discs}

Let $\rho_0$ with $0< \rho_0 \leq \varepsilon_1/4$ and for every
$\rho$ with $0\leq \rho \leq \rho_0$, consider an analytic disc
\def\theequation{9.1}\begin{equation}
A_\rho(\zeta)=(Z_\rho(\zeta), W_\rho(\zeta))=
(X_\rho(\zeta)+iY_\rho(\zeta), 
U_\rho(\zeta)+iV_\rho(\zeta)),
\end{equation}
where $Z_\rho(\zeta)$ is given by
\def\theequation{9.2}\begin{equation}
Z_\rho(\zeta)= (\rho(1-\zeta), 0, \dots, 0).
\end{equation}
The disc $A_\rho$ should be attached to $M\cap\Delta_n(\varepsilon_1)$
and should satisfy $A_\rho(1)= p_1=0$ (remind that $p_1$ is the origin
in the chosen coordinates). A necessary and sufficient condition is
that $U_\rho$ satisfies the so-called Bishop functional equation
\def\theequation{9.3}\begin{equation}
U_\rho(\zeta)=-[T_1(\varphi(X_\rho(\cdot), 
Y_\rho(\cdot), U_\rho(\cdot)))](\zeta), 
\end{equation}
for all $\zeta\in \partial \Delta$.  Here, $T_1$ denotes the harmonic
conjugate operator (Hilbert transform on $\partial \Delta$) normalized
at $\zeta = 1$, namely satisfying $T_1 u(1) = 0$ for every $u \in
\mathcal{ C}^\infty (\partial \Delta, \R^d)$.  By~\cite{ ber}, the
solution exists, is unique and yields a family of analytic discs
$A_\rho(\zeta)$ which is smooth
(and in fact real analytic) with respect to $\rho$ and
$\zeta$. Of course, for $\rho=0$, the disc $A_0$ is the constant disc
which maps $\overline{\Delta}$ to $p_1$.

\bigskip
\begin{center}
\begin{picture}(0,0)%
\includegraphics{disc-wall.pstex}%
\end{picture}%
\setlength{\unitlength}{3947sp}%
\begingroup\makeatletter\ifx\SetFigFont\undefined%
\gdef\SetFigFont#1#2#3#4#5{%
  \reset@font\fontsize{#1}{#2pt}%
  \fontfamily{#3}\fontseries{#4}\fontshape{#5}%
  \selectfont}%
\fi\endgroup%
\begin{picture}(5799,3729)(1189,-3103)
\put(3797, 25){\makebox(0,0)[lb]{\smash{{\SetFigFont{8}{9.6}{\familydefault}{\mddefault}{\updefault}{\color[rgb]{0,0,0}$M_1$}%
}}}}
\put(6113,-297){\makebox(0,0)[lb]{\smash{{\SetFigFont{8}{9.6}{\familydefault}{\mddefault}{\updefault}{\color[rgb]{0,0,0}$M$}%
}}}}
\put(3646,348){\makebox(0,0)[lb]{\smash{{\SetFigFont{8}{9.6}{\familydefault}{\mddefault}{\updefault}{\color[rgb]{0,0,0}$z_2,\ldots,z_m,u$}%
}}}}
\put(4890,-462){\makebox(0,0)[lb]{\smash{{\SetFigFont{8}{9.6}{\familydefault}{\mddefault}{\updefault}{\color[rgb]{0,0,0}$y_1$}%
}}}}
\put(6548,-1302){\makebox(0,0)[lb]{\smash{{\SetFigFont{8}{9.6}{\familydefault}{\mddefault}{\updefault}{\color[rgb]{0,0,0}$x_1$}%
}}}}
\put(2425,-1448){\makebox(0,0)[lb]{\smash{{\SetFigFont{7}{8.4}{\familydefault}{\mddefault}{\updefault}{\color[rgb]{0,0,0}$E$}%
}}}}
\put(3242,-267){\makebox(0,0)[lb]{\smash{{\SetFigFont{8}{9.6}{\familydefault}{\mddefault}{\updefault}{\color[rgb]{0,0,0}$M_1^-$}%
}}}}
\put(4179,-267){\makebox(0,0)[lb]{\smash{{\SetFigFont{8}{9.6}{\familydefault}{\mddefault}{\updefault}{\color[rgb]{0,0,0}$M_1^+$}%
}}}}
\put(3609,-1505){\makebox(0,0)[lb]{\smash{{\SetFigFont{8}{9.6}{\familydefault}{\mddefault}{\updefault}{\color[rgb]{0,0,0}$p_1=A_{\rho_0}(1)$}%
}}}}
\put(4647,-1120){\makebox(0,0)[lb]{\smash{{\SetFigFont{8}{9.6}{\familydefault}{\mddefault}{\updefault}{\color[rgb]{0,0,0}$A_{\rho_0}(\partial\Delta)$}%
}}}}
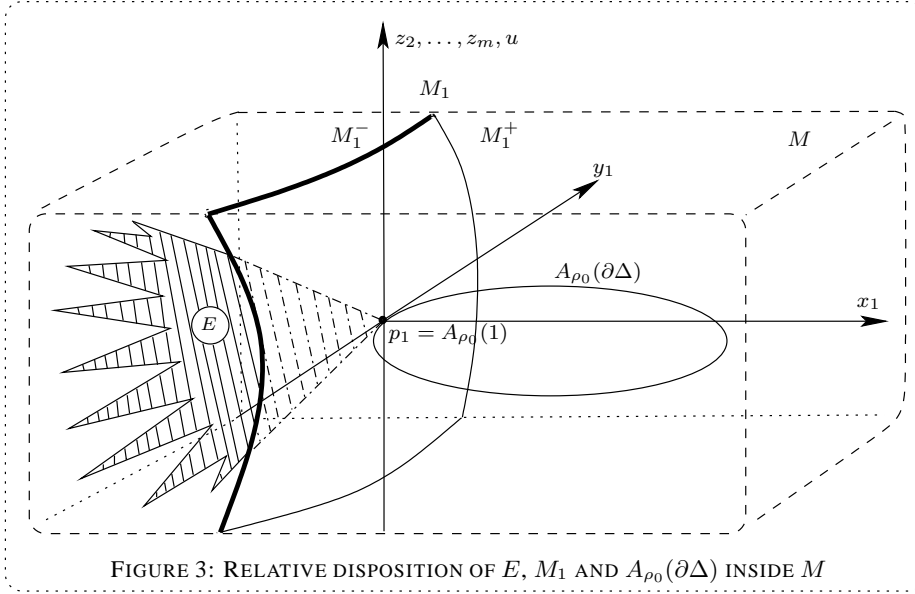
\put(1859,-3001){\makebox(0,0)[lb]{\smash{{\SetFigFont{9}{10.8}{\familydefault}{\mddefault}{\updefault}{\color[rgb]{0,0,0}{\sc Figure~3: Relative disposition of $E$, $M_1$ and $A_{\rho_0}(\partial \Delta)$ inside $M$}}%
}}}}
\end{picture}%

\end{center}
\bigskip

The following elementary lemma shows that the boundary
$A_{\rho_0}(\partial\Delta)$ of the disc $A_{\rho_0}$ meets the wall
$M_1$ only at $p_1$, as shown in {\sc Figure~3} just above.
The notion of analytic isotopy is appropriate
to apply the continuity principle ({\em cf.}~\cite{ mp1}, 
\cite{ mp2}).

\def\thelemma{9.4}\begin{lemma}
There exists $\rho_0$ with $0 < \rho_0 < \varepsilon_1/4$ such that 
the following two properties are satisfied:
\begin{itemize}
\item[{\bf (1)}]
for every $\rho$ with $0< \rho \leq \rho_0$, the mapping
$A_\rho: \overline{\Delta} \to \Delta_n(\varepsilon_1)$ is
an embedding; moreover, each $A_\rho$ is analytically 
isotopic to the point $p_1$;
\item[{\bf (2)}]
for every $\rho$ with $0 < \rho \leq \rho_0$:
\def\theequation{9.5}\begin{equation}
A_\rho(\partial \Delta) \backslash \{1\}) \subset 
M_1^+.
\end{equation}
\end{itemize}
\end{lemma}

\proof
The fact that $A_\rho$ is an embedding for $\rho>0$ is 
obvious, since this is the case for $Z_\rho$. Then the analytic
isotopy is obtained by letting $\rho$ decrease to $0$.
This proves part {\bf (1)}.

For part {\bf (2)}, put $\zeta= re^{i\theta}$, where 
$\vert \theta \vert \leq \pi$ and compute
\def\theequation{9.6}\begin{equation}
\left\{
\aligned
X_{1,\rho}(\zeta)
& \
={\rho(1-\zeta)+\rho(1-\bar \zeta)\over 2}
=
{\rho\over 2} \, 
\left\vert 1-e^{i\theta} \right\vert^2, \\
Y_{1,\rho}(\zeta)
& \
=
{\rho(1-\zeta)-\rho(1-\bar \zeta)\over 2i}
=
-\rho \, \sin \theta.
\endaligned\right.
\end{equation}
Since the real analytic
solution $U_\rho(\zeta)$ vanishes identically for $\rho=0$
and vanishes at $\zeta=1$, there exists a constant $C>0$ such that
\def\theequation{9.7}\begin{equation}
\left\vert U_\rho(\zeta) \right\vert \leq
C \, \rho \, \vert 1 -\zeta \vert,
\end{equation}
for all $\rho\leq \varepsilon_1/4$. One then deduces from the
elementary inequalities ${2 \vert \theta \vert \over \pi}\leq 
\vert 1 -e^{i\theta} \vert \leq \vert \theta \vert$ and 
$\vert \sin \theta \vert \leq \vert \theta \vert$ that
\def\theequation{9.8}\begin{equation}
X_{1,\rho}(e^{i\theta}) \geq {\rho \over 2 \pi^2}\, 
\theta^2, \ \ \ \ \ 
\left\vert Y_{1,\rho}(e^{i\theta}) 
\right\vert^2 \leq \rho^2 \, \theta^2, \ \ \ \ \ 
\left\vert U_\rho(e^{i\theta}) \right\vert^2 \leq C^2 \, \rho^2 \, \theta^2.
\end{equation}
Recall that $Z_{\sharp, \rho}(re^{i\theta})\equiv 0$.
Hence it suffices to choose 
\def\theequation{9.9}\begin{equation}
\rho_0 < {1\over 2\pi^2 (1+C^2)}
\end{equation}
in order to insure that 
\def\theequation{9.10}\begin{equation}
X_{1,\rho}(e^{i\theta}) >
(Y_{1,\rho}(e^{i\theta}))^2+ 
\left\vert U_\rho(e^{i\theta}) \right\vert^2,
\end{equation}
for all $\theta\neq 0$, which completes the proof.
\endproof

\section*{\S10.~Meromorphic 
extension} 

The goal is to construct a meromorphic 
extension of $H_1/ H_2$ to a neighborhood of $p_1$.
Let $\Omega \subset \mathcal{ V}_{\varepsilon_1}$ be an open
neighborhood of the point $A_{\rho_0} (-1)$.  Thanks to Lemma~2.7
in~\cite{ mp1} (a slight modification of the geometric constructions
in~\cite{ tu2}), it is possible to include the disc $A_{\rho_0}
(\zeta)$ in a {\sl regular} 
(in the sense of Definition~1.8 in~\cite{ mp1};
{\it see} also p.~493 of~\cite{ mp2}) family of analytic discs
$A_{\rho_0, s, v}(\zeta)$, where $s\in \R^{2m+d-1}$ satisfies $\vert s
\vert < s_0$ for some $s_0>0$, where $v\in \R^{d-1}$ satisfies $\vert
v \vert < v_0$ for some $v_0>0$ and where $A_{\rho_0, s,v}(\partial
\Delta)$ is contained in
\def\theequation{10.1}\begin{equation}
(M\cap \Delta_n(\varepsilon_1)) \cup 
(\mathcal{W}_{\varepsilon_1} \cap \Omega).
\end{equation}
Essentially, the disc is deformed near $A_{\rho_0}(-1)$ in order that
its direction at exit at $A_{\rho_0}(1)= p_1$ covers an open cone at
$p_1$, by means\ of the parameter $v$, as in~\cite{ tu2}, ~\cite{
mp1}, \cite{ mp2}. Then the parameter $s$ achieves translation along
$M$. By an application of the continuity principle (where property
{\bf (1)} of Lemma 9.4 is needed), there exist $\theta_0$ with
$\theta_0>0$ and $r_0$ with $r_0>0$, $1-r_0< 1$ such that the
following set covered by pieces of analytic discs
\def\theequation{10.2}\begin{equation}
\mathcal{W}_4:=\left\{A_{\rho_0,s,v}(r e^{i\theta}) :\, 
\vert s \vert < s_0, \, 
\vert v \vert < v_0, \, 
\vert \theta \vert < \theta_0, \, 
1-r_0 < r < 1\right\}
\end{equation}
is a (curved) local wedge of edge $M$ at $p_1$ ({\it see} {\sc
Figure~4} just below) to which meromorphic functions in $\widetilde{
\mathcal{ W}}_{ \varepsilon_1} \cup \mathcal{ V}_{ \varepsilon_1}$
extend meromorphically.

\bigskip
\begin{center}
\begin{picture}(0,0)%
\includegraphics{edge-of-the-wedge.pstex}%
\end{picture}%
\setlength{\unitlength}{3947sp}%
\begingroup\makeatletter\ifx\SetFigFont\undefined%
\gdef\SetFigFont#1#2#3#4#5{%
  \reset@font\fontsize{#1}{#2pt}%
  \fontfamily{#3}\fontseries{#4}\fontshape{#5}%
  \selectfont}%
\fi\endgroup%
\begin{picture}(5799,3624)(1189,-3973)
\put(4211,-3221){\makebox(0,0)[lb]{\smash{{\SetFigFont{9}{10.8}{\familydefault}{\mddefault}{\updefault}{\color[rgb]{0,0,0}$A_{\rho_0,s,v}(\Delta)$}%
}}}}
\put(5522,-2155){\makebox(0,0)[lb]{\smash{{\SetFigFont{9}{10.8}{\familydefault}{\mddefault}{\updefault}{\color[rgb]{0,0,0}$\mathcal{V}_{\varepsilon_1}$}%
}}}}
\put(3622,-1562){\makebox(0,0)[lb]{\smash{{\SetFigFont{9}{10.8}{\familydefault}{\mddefault}{\updefault}{\color[rgb]{0,0,0}$E$}%
}}}}
\put(5945,-1986){\makebox(0,0)[lb]{\smash{{\SetFigFont{6}{7.2}{\familydefault}{\mddefault}{\updefault}{\color[rgb]{0,0,0}$M$}%
}}}}
\put(4039,-1497){\makebox(0,0)[lb]{\smash{{\SetFigFont{9}{10.8}{\familydefault}{\mddefault}{\updefault}{\color[rgb]{0,0,0}$p_1$}%
}}}}
\put(2832,-1643){\makebox(0,0)[lb]{\smash{{\SetFigFont{9}{10.8}{\familydefault}{\mddefault}{\updefault}{\color[rgb]{0,0,0}$M_1^-$}%
}}}}
\put(4755,-1602){\makebox(0,0)[lb]{\smash{{\SetFigFont{7}{8.4}{\familydefault}{\mddefault}{\updefault}{\color[rgb]{0,0,0}$M_1^+$}%
}}}}
\put(5206,-2446){\makebox(0,0)[lb]{\smash{{\SetFigFont{9}{10.8}{\familydefault}{\mddefault}{\updefault}{\color[rgb]{0,0,0}$A_{\rho_0}(\Delta)$}%
}}}}
\put(5696,-1838){\makebox(0,0)[lb]{\smash{{\SetFigFont{5}{6.0}{\familydefault}{\mddefault}{\updefault}{\color[rgb]{0,0,0}$A_{\rho_0}$}%
}}}}
\put(5891,-1838){\makebox(0,0)[lb]{\smash{{\SetFigFont{5}{6.0}{\familydefault}{\mddefault}{\updefault}{\color[rgb]{0,0,0}$(-1)$}%
}}}}
\put(2823,-621){\makebox(0,0)[lb]{\smash{{\SetFigFont{9}{10.8}{\familydefault}{\mddefault}{\updefault}{\color[rgb]{0,0,0}$\Delta_n(\varepsilon_1)$}%
}}}}
\put(2401,-3292){\makebox(0,0)[lb]{\smash{{\SetFigFont{9}{10.8}{\familydefault}{\mddefault}{\updefault}{\color[rgb]{0,0,0}$\widetilde{\mathcal{W}}_{\varepsilon_1}$}%
}}}}
\put(3611,-1166){\makebox(0,0)[lb]{\smash{{\SetFigFont{9}{10.8}{\familydefault}{\mddefault}{\updefault}{\color[rgb]{0,0,0}$\widetilde{\mathcal{W}}_4$}%
}}}}
\put(3762,-2112){\makebox(0,0)[lb]{\smash{{\SetFigFont{9}{10.8}{\familydefault}{\mddefault}{\updefault}{\color[rgb]{0,0,0}$\mathcal{W}_4$}%
}}}}
\put(2026,-1994){\makebox(0,0)[lb]{\smash{{\SetFigFont{9}{10.8}{\familydefault}{\mddefault}{\updefault}{\color[rgb]{0,0,0}$M$}%
}}}}
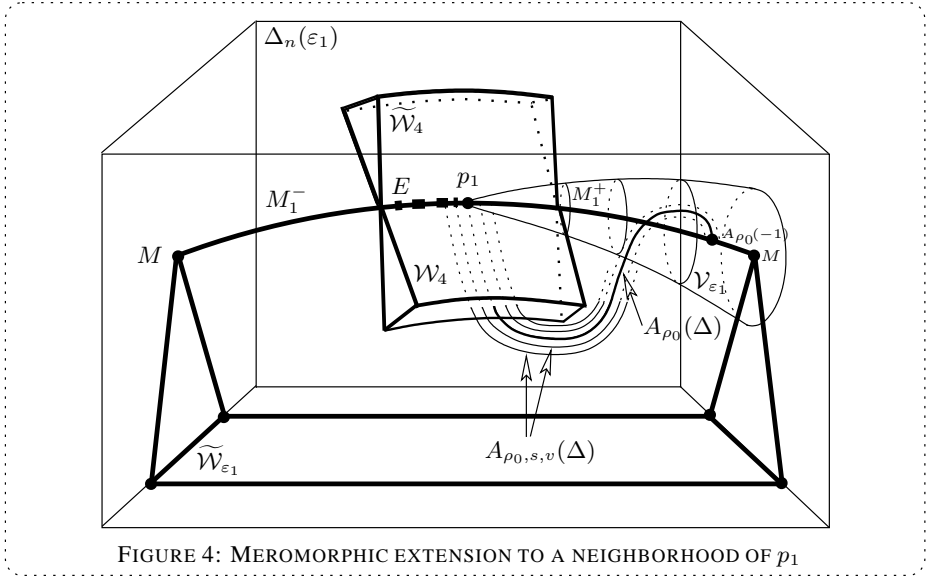
\put(1903,-3876){\makebox(0,0)[lb]{\smash{{\SetFigFont{9}{10.8}{\familydefault}{\mddefault}{\updefault}{\color[rgb]{0,0,0}{\sc Figure~4: Meromorphic extension to a neighborhood of $p_1$}}%
}}}}
\end{picture}%

\end{center}
\bigskip

\def\thelemma{10.3}\begin{lemma}
In the preceding situation, three extension results hold:
\begin{itemize}
\item[{\bf (1)}]
for $j=1,2$, the $\C^{\nu_j}$-valued $\mathcal{C}^\infty$-smooth CR
function $b_j$ extends holomorphically to $\mathcal{W}_4$;
\item[{\bf (2)}]
the CR quotient $H_1/H_2$ extends meromorphically to a
symmetric wedge $\widetilde{\mathcal{W
}}_4$ contained in $s(\mathcal{W}_4)$
which is only slightly smaller;
\item[{\bf (3)}]
the CR quotient $H_1/H_2$ extends meromorphically to 
$\mathcal{W}_4$.
\end{itemize}
\end{lemma}

\proof
As a preliminary, define translations of geometric objects in the
normal directions $T_{p_1}\C^n/T_{p_1}M$ as follows. 
If a unitary vector
$\upsilon_1\in T_{p_1}\C^n$ with zero $z$-component and zero
$u$-component is given, namely $\upsilon_1$ is of the form $(0,
iv_1)\in\C^m\times \C^d$ with $\vert v_1 \vert =1$, 
then for every $\eta$
very small in comparison with $\varepsilon_1$, define the
translation
\def\theequation{10.4}\begin{equation}
(M\cap \Delta_n(\varepsilon_1))+\eta\, \upsilon_1.
\end{equation}
As $\mathcal{ W}_{ \varepsilon_1}$ and $\widetilde{ \mathcal{W}}_{
\varepsilon_1}$ are (approximatively) symmetric to each other, there
exists a unitary vector $\upsilon_1=(0,iv_1)\in T_{p_1} 
\C^n$ such that 
\def\theequation{10.5}\begin{equation}
\left\{
\aligned
(M\cap \Delta_n(\varepsilon_1))+\eta \, \upsilon_1 \subset 
\mathcal{W}_{\varepsilon_1} \ \ \ \ \ \ {\rm for} \ \ \eta>0,\\
(M\cap \Delta_n(\varepsilon_1))+\eta \, \upsilon_1 \subset 
\widetilde{\mathcal{W}}_{\varepsilon_1} \ \ \ \ \ \ {\rm for} \ \ \eta<0.
\endaligned\right.
\end{equation}
In other words, $\upsilon_1 \in T_{ p_1} \mathcal{ W}_{
\varepsilon_1}$ and $-\upsilon_1 \in T_{p_1} \widetilde{ \mathcal{
W}}_{\varepsilon_1}$. Translate also the analytic discs, which yields
$A_{\rho_0, s, v}(\zeta)+\eta\, \upsilon_1$.

Prove now {\bf (1)} of Lemma~10.3.  By construction, the two $\mathcal{
C }^\infty$-smooth CR functions $b_1$ and $b_2$ extend holomorphically
to $\mathcal{ W}_{ \varepsilon_1}$. Since for every $\eta>0$, the
discs $A_{\rho_0,s,v}(\zeta )+\eta\, \upsilon_1$ have their boundaries
contained in $\mathcal{ W }_{ \varepsilon_1}$ and are analytically
isotopic to a point in $\mathcal{ W}_{ \varepsilon_1}$, it follows
from the continuity principle that $b_1$ and $b_2$ extend
holomorphically to the wedge $\mathcal{W}_4+\eta\, \upsilon_1$. By
letting $\eta$ tend to zero, it follows that $b_1$ and $b_2$ extend
holomorphically to $\mathcal{W}_4$.

Next, prove {\bf (2)} of Lemma~10.3.  Since $b_1$ and $b_2$
extend holomorphically to $\mathcal{W}_4$, by reasoning as in the
beginning of the proof of Main Lemma~5.2, it follows that $H_1/H_2$
extends meromorphically to the symmetric wedge
$\widetilde{\mathcal{W}}_4$, which is
contained in $s(\mathcal{W}_4)$, but only
slightly smaller. Importantly, there exists a unitary vector
$\upsilon_4 \in T_{p_1}\C^n$ with coordinates of the form $(0,iv_4)\in
\C^m\times \C^d$ such that $\upsilon_4 \in T_{p_1}\mathcal{W}_4$ and
$-\upsilon_4 \in T_{p_1} \widetilde{\mathcal{W}}_4$
({\it see} again {\sc Figure~4}).

Finally, prove {\bf (3)} of Lemma~10.3. In the domain
$\widetilde{ \mathcal{W}}_{\varepsilon_1} \cup \mathcal{ V}_{
\varepsilon_1}$, the meromorphic extension of the CR quotient $H_1 /
H_2$, can be represented as a quotient $P_1/P_2$.
Since for every $\eta<0$, the discs $A_{ \rho_0,s,v}(\zeta)+ \eta\,
\upsilon_1$ have their boundaries contained in $\widetilde{ \mathcal{
W}}_{ \varepsilon_1} \cup \mathcal{V}_{ \varepsilon_1}$ and are
analytically isotopic to a point in $\widetilde{\mathcal{
W}}_{\varepsilon_1} \cup \mathcal{V}_{ \varepsilon_1}$, it follows
from the continuity principle that $P_1$ and $P_2$
extend holomorphically to $\mathcal{ W}_4+ \eta\, \upsilon_1$. By
letting $\eta$ tend to zero, one deduces that $P_1$ and $P_2$ extend
holomorphically to $\mathcal{W}_4$.  Hence the CR quotient $H_1 / H_2$
extends meromorphically to $\mathcal{W}_4$.

The proof of Lemma~10.3 is complete.
\endproof

The proof of Main Lemma~7.1 is almost achieved. By
Lemma~10.3, the CR quotient $H_1/H_2$ extends 
meromorphically to the union
\def\theequation{10.6}\begin{equation}
\mathcal{W}_4 \cup \widetilde{
\mathcal{W}}_4 \cup \mathcal{V}_{
\varepsilon_1}.
\end{equation}
This union is connected ({\it cf.}~{\sc Figure~5} just below.  Denote
again by $P_1/P_2$ this meromorphic extension, where $P_1$ and $P_2$
are holomorphic in the domain $\mathcal{W}_4 \cup \widetilde{
\mathcal{W}}_4 \cup \mathcal{V}_{ \varepsilon_1}$.

\bigskip
\begin{center}
\begin{picture}(0,0)%
\includegraphics{deformation.pstex}%
\end{picture}%
\setlength{\unitlength}{3947sp}%
\begingroup\makeatletter\ifx\SetFigFont\undefined%
\gdef\SetFigFont#1#2#3#4#5{%
  \reset@font\fontsize{#1}{#2pt}%
  \fontfamily{#3}\fontseries{#4}\fontshape{#5}%
  \selectfont}%
\fi\endgroup%
\begin{picture}(5799,2551)(1192,-3234)
\put(5522,-2155){\makebox(0,0)[lb]{\smash{{\SetFigFont{9}{10.8}{\familydefault}{\mddefault}{\updefault}{\color[rgb]{0,0,0}$\mathcal{V}_{\varepsilon_1}$}%
}}}}
\put(3622,-1562){\makebox(0,0)[lb]{\smash{{\SetFigFont{9}{10.8}{\familydefault}{\mddefault}{\updefault}{\color[rgb]{0,0,0}$E$}%
}}}}
\put(5945,-1986){\makebox(0,0)[lb]{\smash{{\SetFigFont{6}{7.2}{\familydefault}{\mddefault}{\updefault}{\color[rgb]{0,0,0}$M$}%
}}}}
\put(3611,-1166){\makebox(0,0)[lb]{\smash{{\SetFigFont{9}{10.8}{\familydefault}{\mddefault}{\updefault}{\color[rgb]{0,0,0}$\widetilde{\mathcal{W}}_4$}%
}}}}
\put(2026,-1994){\makebox(0,0)[lb]{\smash{{\SetFigFont{9}{10.8}{\familydefault}{\mddefault}{\updefault}{\color[rgb]{0,0,0}$M$}%
}}}}
\put(3845,-2874){\makebox(0,0)[lb]{\smash{{\SetFigFont{9}{10.8}{\familydefault}{\mddefault}{\updefault}{\color[rgb]{0,0,0}$\mathcal{W}_4^d$}%
}}}}
\put(3832,-1469){\makebox(0,0)[lb]{\smash{{\SetFigFont{9}{10.8}{\familydefault}{\mddefault}{\updefault}{\color[rgb]{0,0,0}$M^d$}%
}}}}
\put(4260,-1194){\makebox(0,0)[lb]{\smash{{\SetFigFont{9}{10.8}{\familydefault}{\mddefault}{\updefault}{\color[rgb]{0,0,0}$p_1$}%
}}}}
\put(5206,-2446){\makebox(0,0)[lb]{\smash{{\SetFigFont{9}{10.8}{\familydefault}{\mddefault}{\updefault}{\color[rgb]{0,0,0}$A_{\rho_0}^d(\Delta)$}%
}}}}
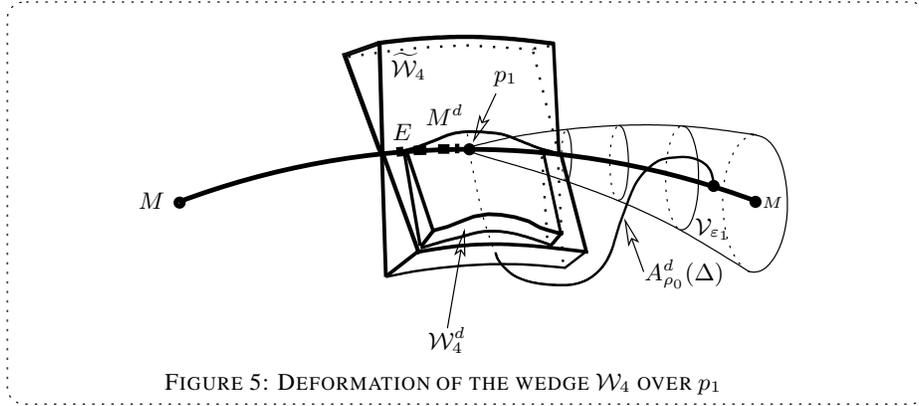
\put(2190,-3130){\makebox(0,0)[lb]{\smash{{\SetFigFont{9}{10.8}{\familydefault}{\mddefault}{\updefault}{\color[rgb]{0,0,0}{\sc Figure~5:  Deformation of the wedge $\mathcal{W}_4$ over $p_1$}}%
}}}}
\end{picture}%

\end{center}
\bigskip

Introduce a one parameter family of smooth deformations $M^d$, $d\geq
0$, of $M$ localized in a neighborhood of $p_1$, with $M^0=M$, by
pushing $M$ near $p_1$ inside $\widetilde{ \mathcal{ W}}_4$. For this,
it suffices to
use a cut-off function $\chi(x, y,u)$ with support in a
neighborhood of the origin, and for $d\in\R$, $d\geq 0$, to define
$M^d$ by the vectorial equation
\def\theequation{10.7}\begin{equation}
v=\varphi(x,y,u)+d\, \chi(x,y,u) \, \upsilon_4,
\end{equation} 
where $\upsilon_4 \in T_{p_1}\mathcal{W}_4$ and $-\upsilon_4 \in
T_{p_1} \widetilde{\mathcal{W}}_4$.

Since the resolution of Bishop's equation is stable under
perturbation, there exists a smoothly deformed family $A_{\rho_0, s,v
}^d( \zeta)$ of analytic discs with $A_{\rho_0, s,v }^0( \zeta)=A_{
\rho_0,s,v}(\zeta)$ and a deformed wedge
\def\theequation{10.8}\begin{equation}
\mathcal{W}_4^d:=\{
A_{\rho_0,s,v}^d(r e^{i\theta}) : \, 
\vert s \vert < s_0, \, 
\vert v \vert < v_0, \, 
\vert \theta \vert < \theta_0, \, 
1-r_0 < r < 1\}.
\end{equation}
Since $\mathcal{W}_4$ and $\widetilde{ \mathcal{ W}}_4$ have opposite
directions, it is clear that one can insure that $\mathcal{ W}_4^d$
contains the point $p_1$.

For $d>0$, let the parameter $\rho$ vary in the interval $[0, \rho_0]$
to deduce that the discs $A_{\rho_0, s,v}^d$ are analytically isotopic
to a point in $\mathcal{ W}_4\cup \widetilde{ \mathcal{ W}}_4 \cup
\mathcal{ V}_{ \varepsilon_1}$. By a further application of the
continuity principle, it follows that $P_1$ and $P_2$ extend
holomorphically to $\mathcal{ W}_4^d$, hence to a neighborhood of
$p_1$.  In conclusion, the CR quotient $H_1/H_2$
extends meromorphically to a neighborhood of $p_1$.

The proof of Main Lemma~7.1 is complete. 
\endproof

In conclusion, the proof of Theorem~2.9 is complete.

\section*{\S11~Proof of Lemma~2.12} 

Establish first the second sentence. Assume that the source $M$ is a
hypersurface and that $h$ is essentially finite at $p_0$.  Proceeding
by contradiction, suppose that $M$ is not minimal at $p_0$.
Equivalently, the CR orbit of $p_0$ is an $(n-1)$-dimensional complex
hypersurface passing through $p_0$, which coincides in fact with the
Segre variety $S_{\bar p_0}$.  Choose holomorphic coordinates $(z,w)
\in\C^{n-1}\times \C$ vanishing at $p_0$ in which $S_{\bar p_0} = \{
w= 0\}$ and represent $M$ by a single complex equation $\bar
w=\Theta(\bar z,z,w)$, where $S_{\bar p_0}=\{(z,0)\}$, whence
$\Theta(\bar z,z,0)\equiv 0$. It follows that the usual basis of
$(0,1)$ vector fields tangent to $M$ satisfies
$\overline{L}_k=\partial /\partial \bar z_k+ {\rm O}(w)\, \partial /
\partial \bar w$.

In the fundamental equations
\def\theequation{11.1}\begin{equation}
r_{j'}'(h(t), \overline{h(t)})=0,
\end{equation}
for $j'=1,\dots,d'$, specify $t=(z,0)\in M$, which yields:
\def\theequation{11.2}\begin{equation}
r_{j'}'(h(z,0),\bar h(\bar z,0))\equiv 0, \ \ \ \ \ 
j'=1,\dots,d'.
\end{equation}
Since $h$ is CR and $\{(z,0)\}$ is contained in $M$, 
the mapping $z \mapsto h(z,0)$ is holomorphic,
hence it
is justified to write $\bar h(\bar z,0)$ instead of 
$\overline{h(z,0)}$.

By assumption, $h$ is essentialy finite 
at $p_0$. Hence apply the
derivations $\overline{L}^\beta$, for $\beta\in\N^{n-1}$,
to~\thetag{11.1}, which amounts to applying
the derivations $\partial_{\bar z}^\beta$ to~\thetag{11.2}, 
which yields expressions of the form
\def\theequation{11.3}\begin{equation}
S_{j',\beta}'\left( 
J_{\bar z}^{\vert \beta \vert} \bar h (\bar z,0) :\,
h(z,0) \right)\equiv 0.
\end{equation}
Exactly as in \S3, it follows 
from the essential 
finiteness assumption
that there exist Weierstrass polynomials such that
\def\theequation{11.4}\begin{equation}
(h_{i'}(z,0))^{N_{i'}'}+
\sum_{1\leq I' \leq N_{i'}'}\,
H_{i',I'} \left(
J_{\bar z}^{\ell_0} \bar h(\bar z,0)\right)\, 
(h_{i'}(z,0))^{N_{i'}'-I'}= 0,
\end{equation}
Putting $\bar z=0$, it follows that $h_{i'}(z,0)$ is a constant for
$i'=1,\dots, n'$, hence vanishes identically. However, if $h(z,0)$
vanishes identically, it is clearly impossible for $h$ to be
essentially finite at $p_0$, since differentiation with respect to
$\bar z$ of $r_{j'}'(t',\bar h(\bar z,0))$ gives nothing else than the
constant zero, so $\V_0'$ is defined by $\{ w' =0\}$, hence is
positive-dimensional.  This establishes the
second sentence of Lemma~2.12.

To establish the first sentence, remind again that 
if $\dim \, \mathcal{O}_{CR}(M',p_0')=2m'$, then the CR orbit of
$p_0'$ is an $m'$-dimensional complex manifold passing through $p_0'$,
which coincides in fact with the Segre variety $S_{\bar
p_0'}'$. Assume that the coordinates $(z',w')\in\C^{m'}\times \C^{d'}$
are such that $S_{\bar p_0'}' =\{ w'=0\}$. Then the complex defining
equations of $M'$ are of the form $\bar w_{j'}'= \Theta_{j'}'(\bar
z',z',w')$, where $\Theta_{j'}'(\bar z',z',0)\equiv 0$ for
$j'=1,\dots,d'$. However, this clearly contradicts the two classical
characterizations of essential finiteness (\cite{ ber}).

In general, there exist
coordinates $(z',w_1',w_2')\in\C^{m'}\times \C^{d_1'}\times
\C^{d_2'}$ centered at $p_0'$, 
where $d_2'$ is the holomorphic codimension of the
intrinsic complexification $(\mathcal{O}_{CR}(M',p_0'))^{i_c}$ and
where $d_1'=d'-d_2'$,
such that $M'$ is represented by
\def\theequation{11.5}\begin{equation}
\left\{
\aligned
\bar w_{1,j_1'}=\Theta_{1,j_1'}'(\bar z, z', w_1', w_2'), 
\ \ \ \ \ \ j_1'=1,\dots,d_1',\\
\bar w_{2,j_2'}=\Theta_{2,j_2}'(\bar z', z', w_1', w_2'),
\ \ \ \ \ \ j_2'=1,\dots,d_2',
\endaligned\right.
\end{equation}
where $\Theta_{2',j_2'}'(\bar z',z',w_1',0)\equiv 0$ for
$j_2'=1,\dots,d_2'$.  
Without loss of generality, assume that the
coordinates are normal.
By the characterization of essential finiteness in normal coordinates,
the ideal in $\C \{ z'\}$ generated by the partial derivatives
$\partial_{\bar z'}^{\beta'} \Theta_{1,j_1'}' (z',0, 0, 0)$, 
$j_1 = 1,\dots, d_1'$, and the
partial derivatives $\partial_{\bar z'}^{\beta'}
\Theta_{2,j_2'}' (z',0, 0, 0)$, $j_2' = 1, \dots, d_2'$, 
should be of finite
codimension, where 
$\beta'$ runs in $\N^{ m'}$. 
However, the second collection vanishes identically.
Thus, the ideal in $\C \{ z'\}$ generated by the partial
derivatives $\partial_{\bar z'}^{\beta'} \Theta_{1,j_1'}' (z',0, 0,
0)$ is of finite codimension. This shows that the CR orbit $M'\cap
\{w_2'=0\}$ is essentially finite at the origin.  The proof of
Lemma~2.12 is complete.

\vfill
\end{document}